\RequirePackage{fix-cm}
\documentclass[twocolumn]{svjour3}          

\usepackage{graphicx}
\usepackage{amsmath}
\usepackage{amsthm}
\usepackage{cite}
\usepackage{epstopdf}
\usepackage{caption}
\usepackage{epstopdf}
\usepackage{amsfonts}
\usepackage{quotes}
\usepackage{fancyhdr}
\usepackage{subfigure}
\usepackage{amssymb}
\usepackage{color}
\usepackage{algorithmic}
\usepackage[margin=1in]{geometry}
\usepackage{amssymb}
\usepackage{framed}
\usepackage{varwidth}
\usepackage{mathtools}

\DeclareMathOperator*{\argmax}{arg\,max}
\DeclareMathOperator*{\argmin}{arg\,min}

\newcommand{\envert}[1]{\left\lvert#1\right\rvert}
\let\abs=\envert

\newcommand{\enVert}[1]{\left\lVert#1\right\rVert}
\let\norm=\enVert

\newcommand{\inProd}[1]{\left\langle#1 \right\rangle}

\DeclareMathOperator{\diver}{div}

\newcommand{\bitem}{\begin{itemize}}
\newcommand{\eitem}{\end{itemize}}

\newcommand{\bpm}{\begin{pmatrix}}
\newcommand{\epm}{\end{pmatrix}}

\newcommand{\bq}{\begin{equation}}
\newcommand{\eq}{\end{equation}}

\usepackage{etoolbox}
\patchcmd{\theorem}{Theorem}{Teorema}{}{}

\usepackage{color}

\begin{document}

\title{Convex variational methods on graphs for multiclass segmentation of high-dimensional data and point clouds}

\author{Egil Bae \and Ekaterina Merkurjev}

\institute{Egil Bae  \at
Norwegian Defence Research Establishment\\ P.O. Box 25, NO-2027 Kjeller, Norway \\ \email{Egil.Bae@ffi.no}
           \and
          Ekaterina Merkurjev \at
Michigan State University \\
                220 Trowbridge Rd, East Lansing, MI 48824 \\
          \email{kmerkurev@math.msu.edu } 
}

\maketitle

\begin{abstract}

Graph-based variational methods have recently shown to be highly competitive for various classification problems of high-dimensional data, but are inherently difficult to handle from an optimization perspective. This paper proposes a convex relaxation for a certain set of graph-based multiclass data segmentation models involving a graph total variation term, region homogeneity terms, supervised information and certain constraints or penalty terms acting on
the class sizes. Particular applications include semi-supervised classification of high-dimensional data and unsupervised segmentation of unstructured 3D point clouds. Theoretical analysis shows that the convex relaxation closely approximates the original NP-hard problems, and these observations are also confirmed experimentally. An efficient duality based algorithm is developed that handles all constraints on the labeling function implicitly. Experiments on semi-supervised classification indicate consistently higher accuracies than related non-convex approaches, and considerably so when the training data are not uniformly distributed among the data set. The accuracies are also highly competitive against a wide range of other established methods on three benchmark datasets. Experiments on 3D point clouds acquired by a LaDAR in outdoor scenes, demonstrate that the scenes can accurately be segmented into object classes such as vegetation, the ground plane and human-made structures.


\end{abstract}


\section{Introduction}\label{Introduction}

The graphical framework has become a popular setting for classification \cite{zhou:bousquet:lal,zhou:scholkopf,wang,belkin,chapelle:scholkopf:zien,zhu} and filtering \cite{Di12,ETT15,TMR09,MBS13,TS16,TSBLB14} of high-dimensional data. 
Some of the best performing classification algorithms are based on solving variational problems on graphs \cite{szlam,bresson:laurent,bertozzi,merkurjev,HB10,BTCS14,bresson_2013,MBBT15,gilboa2,zhou:bousquet:lal,chapelle:scholkopf:zien,AGP15, TEL14, HLE13}. In simple terms, these algorithms attempt to group the data points into classes in such a way that pairs of data points with different class memberships are as dissimilar as possible with respect to a certain feature.
In order to avoid the computational complexity of working with fully connected graphs, approximations, such as those based on spectral graph theory \cite{bertozzi,merkurjev,garcia} or nearest neighbors \cite{ELB08,BTCS14,MBBT15}, are typically employed. For example, \cite{bertozzi} and \cite{merkurjev} employ spectral approaches along with the Nystr\"{o}m extension method \cite{fowlkes} to efficiently calculate the eigendecomposition of a dense graph Laplacian. Works, such as \cite{ELB08,gilboa2,citeulike:10200012,DBLP:journals/siamis/ZhangBBO10,BTCS14,MBBT15}, use the "nearest neighbor" approach to sparsify the graph for computational efficiency. Variational problems on graphs have also become popular for processing
of 3D point clouds \cite{ELB08,LEL12,DEL13,HLE13,LEL14,LLZ13}.  


When the classification task is cast as the minimization of similarity of point pairs with different class membership, extra
information is necessary to avoid the trivial global minimizer of value zero where all points are assigned to the same class. In semi-supervised classification methods, a small set of the data points are given as training data in
advance, and their class memberships are imposed as hard constraints in the optimization problem.
In unsupervised classification methods one typically enforces the sizes of each class not to deviate too far from each other, examples including the normalized cut \cite{shi:malik} and Cheeger ratio cut problems \cite{Ch70}.

Most of the computational methods for semi-supervised and unsupervised classification obtain the solution by computing the local minimizer of a non-convex energy functional. Examples of such algorithms are those based on phase fields \cite{bertozzi} and the MBO scheme \cite{merkurjev,garcia,merkurjev_aml,sunu}. PDEs on graphs for semi-supervised classification also include the Eikonal equation \cite{DEL13} and tug-of-war games related to the infinity-Laplacian equation \cite{ETT15}. Unsupervised problems with class size constraints are inherently the most difficult to handle from an optimization viewpoint, as the convex envelope of the problem has a trivial constant function as a minimizer \cite{shi:malik,szlam,bresson:laurent}. Various ways of simplifying the energy landscape have been proposed \cite{hein,BTCS14,TSBLB14}.
Our recent work \cite{MBBT15} showed that semi-supervised classification problems with two classes could be formulated in a completely convex framework and also presented efficient algorithms that could obtain global minimizers.
%

Image segmentation is a special classification problem where the objective is to assign each pixel to a region. Algorithms based on energy minimization are among the most successful image segmentation methods, and they have historically been divided into \lq region-based' and \lq contour-based'.

Region-based methods attempt to find a partition of the image so that the pixels within each region as a whole are as similar as possible. Additionally, some regularity is imposed on the region boundaries to favor spatial grouping of the pixels. The similarity is usually measured in the statistical sense. In the simplest case, the pixels within each region should be similar to the mean intensity of each region, as proposed in the Chan-Vese \cite{CV01} and Mumford-Shah \cite{MS89} models. Contour-based methods \cite{KWT88,CKS97} instead seek the best suited locations of the region boundaries, typically at locations of large jumps in the image intensities, indicating the interface between two objects.

More recently, it has been shown that the combination of region and contour-based terms in the energy function can give qualitatively very good results \cite{bresson2007fast,gilboa2,DBLP:journals/siamis/JungPC12}, especially when non-local operators are used in the contour terms \cite{gilboa2,DBLP:journals/siamis/JungPC12,DEL13}. There now exists efficient algorithms for solving the resulting optimization problems that can avoid getting stuck in a local minimum, including both combinatorial optimization algorithms \cite{boykov_kolmogorov,boykov:veksler,Kolmogorov04whatenergy} and more recent convex continuous optimization algorithms \cite{bresson2007fast,ZGFN08,Lellmann-et-al-09a,Pock-et-al-iccv09,CCP12,BYT2011,yuan,YBTB2013,BT15}. The latter have shown to be advantageous in several aspects, such as the fact that they require less memory and have a greater potential for parallel implementation of graphics processing units (GPUs), but special care is needed in case of non-local variants of the differential operators (e.g. \cite{DBLP:conf/eccv/RanftlBP14}).

Most of the current data segmentation methods \cite{szlam,bresson:laurent,bertozzi,merkurjev,HB10,BTCS14,bresson_2013,MBBT15} can be viewed as \lq contour-based', since they seek an optimal location of the boundaries of each region. Region-based variational image segmentation models with two classes were generalized to graphs for data segmentation in \cite{lezoray} and for 3D point cloud segmentation in \cite{lezoray,LEL14,TLE15} in a convex framework.
The region terms could be constructed directly from the point geometry and/or be constructed from a color vector defined at the points. 
Concrete examples of the latter were used for experiments on point cloud segmentation.
Region terms have also been proposed in the context of Markov Random Fields for 3D point cloud segmentation \cite{DBLP:conf/cvpr/AnguelovTCKGHN05,MVH08,Triebel06robust3d},
where the weights were learned from training data using associate Markov networks.
The independent preprint \cite{YTS16} proposed to use region terms for multiclass semi-supervised classification in a convex manner, where the region terms were inferred from the supervised points by diffusion.

\subsubsection*{Contributions}

This paper proposes a convex relaxation and an efficient algorithmic optimization framework for a
general set of graph based data classification problems that exhibits non-trivial global minimizers.
It extends the convex approach for semi-supervised classification with two classes given in our previous work \cite{MBBT15} to a much broader range of problems, including multiple classes, novel and more practically useful incorporation of class size information, and a novel unsupervised segmentation model for 3D point clouds acquired by a LaDAR.

The same basic relaxation for semi-supervised classification also appeared in the independent preprint \cite{YTS16}.
The most major distinctions of this work compared to the preprint \cite{YTS16} are: 
we also incorporate class size information in the convex framework;
we give a mathematical and experimental analysis of the close relation between the convex
relaxed and original problems; we propose a different duality based \lq max-flow' inspired algorithm; we incorporate information of the supervised points in a different way; and we consider unsupervised segmentation of 3D point clouds.

The contributions can be summarized more specifically as follows:
\begin{itemize}

\item We specify a general set of classification problems that are suitable for being approximated in a convex manner.
    The general set of problems involves minimization of a multiclass graph cut term together with
supervised constraints, region homogeneity terms and novel constraints or penalty terms acting on the class sizes.
Special cases include semi-supervised classification of high-dimensional data and unsupervised segmentation of 3D point clouds.
\vspace{0.05in}

\item

A convex relaxation is proposed for the general set of problems and its approximation properties are analyzed thoroughly in theory and experiments. This extends the work on multiregion image segmentation \cite{ZGFN08,Lellmann-et-al-09a,BYT2011} to data clustering on graphs, and to cases where there are constraints or penalty terms acting on the class sizes. Since either the introduction of multiple classes or size constraints causes the general problem to become NP-hard,
the relaxation can (probably) not be proven to be exact. Instead, conditions are derived for when an exact global minimizer
can be obtained from a dual solution of the relaxed problem. The strongest conditions are derived in case there are no constraints on the class sizes, but the theoretical results in both cases show that very close approximations are expected. These theoretical results also
agree well with experimental observations. \vspace{0.05in}

\item

The convex relaxed problems are formulated as
equivalent dual problems that are structurally similar to the \lq max-flow' problem over the graph.
This extends our work \cite{MBBT15} to multiple classes and the work on image segmentation proposed in \cite{YBTB10} to data clustering on graphs. We use a conceptually different proof than \cite{MBBT15,YBTB10}, which relates \lq max-flow' with another more direct dual formulation of the problem. Furthermore, it is shown that also the size constraints and penalty
term can be incorporated naturally in the max-flow problem
by modifying the flow conservation condition, such that there should be a constant flow
excess at each node. \vspace{0.05in}

\item

As in our previous work \cite{MBBT15,YBTB10}, an augmented Lagrangian algorithm is developed
based on the new \lq max-flow' dual formulations of the problems. A key advantage compared to
related primal-dual algorithms \cite{DBLP:journals/jmiv/ChambolleP11} in imaging, such as the one considered in the preprint \cite{YTS16}, is that all constraints
on the labeling function are handled implicitly.
This includes constraints on the class sizes,
which are dealt with by separate dual variables
indicating the flow excess at the nodes. Consequently, projections onto the constraint set of the labeling function,
which tend to decrease the accuracy and put strict
restrictions on the step sizes, are avoided.\vspace{0.05in}

\item

We propose an unsupervised segmentation model for unstructured 3D point clouds aquired by a LaDAR within the general framework.
It extends the models of \cite{lezoray,LEL14,TLE15} to multiple classes and gives concrete examples of region terms constructed purely based on geometrical information of the unlabeled points, in order to distinguish classes such as vegetation, the ground plane and human-made structures in an outdoor scene. We also propose a graph total variation term that favors alignment of the region boundaries along "edges" indicated by discontinuities in the normal vectors or the depth coordinate. In contrast to \cite{DBLP:conf/cvpr/AnguelovTCKGHN05,MVH08,Triebel06robust3d,Golovinskiy:2009:SRO}, our model does not rely on any training data.\vspace{0.05in}

\item

Extensive experimental evaluations on \newline semi-supervised classification indicate consistently higher accuracies than related local minimization approaches, and considerably so when the training data are not uniformly distributed among the data set. The accuracies are also highly competitive against a wide range of other established methods on three benchmark datasets. The accuracies can be improved further if an estimate of the approximate class sizes are given in advance. Experiments on 3D point clouds acquired by a LaDAR in outdoor scenes demonstrate that the scenes can accurately be segmented into object classes such as vegetation, the ground plane and regular structures.
The experiments also demonstrate fast and highly accurate convergence of the algorithms, and show that the approximation difference between the convex and original problems vanishes or becomes extremely low in practice.
\end{itemize}

\subsubsection*{Organization}

This paper starts by formulating the general set of problems mathematically in Section \ref{Segmentation as optimization problem over graph}. Section \ref{Binary representation and convex relaxation of multi-classification problems} formulates a convex relaxation of the general problem and analyzes the quality of the relaxation from a dual perspective. Section \ref{Max-flow dual formulations and algorithms} reformulates the dual problem as a \lq max-flow' type of problem and derives an efficient algorithm. Applications to semi-supervised classification of high-dimensional data are presented in Section \ref{Semi-supervised Classification Results}, and applications to segmentation of unstructured 3D point clouds are described in Section \ref{Segmentation of point clouds}, including specific constructions of each term in the general model.  
Section \ref{Applications and experiments} also presents a detailed experimental evaluation for both applications.

\newpage

\section{Data segmentation as energy minimization over a graph}\label{Segmentation as optimization problem over graph}

Assume we are given $N$ data points in $\mathbb{R}^M$. In order to formulate the segmentation of the data points as a minimization problem, the points are first organized in an undirected graph. Each data point is represented by a node in the graph.
The edges in the graph, denoted by $E$, consist of pairs of data points.
Weights $w(x,y)$ on the edges $(x,y) \in E$ measure the similarity between data points $x$ and $y$. A high value of $w(x,y)$ indicates that $x$ and $y$ are similar and a low value indicates that they are dissimilar. A popular choice for the weight function is the Gaussian
\begin{equation}
            w(x,y)= e^{-\frac{d(x,y)^{2}}{\sigma^{2}}},
\label{w1}
\end{equation}
where $d(x,y)$ is the distance, in some sense, between $x$ and $y$. For example, the distance between two 3D points $x$ and $y$ is naturally their Euclidean distance.
In order to reduce the computational burden of working with fully connected graphs, one often only considers the set of edges where $w(x,y)$ is largest. For instance, edges may be constructed between each vertex in $V$ and its $k$ nearest neighbors. More formally, for each $x \in V$, one constructs an edge $(x,y) \in E$ for the $k$ points with the shortest distance $d(x,y)$ to $x$.
 Such a graph can be constructed efficiently by using kd-trees in $O(N k\, \text{log}(N k))$ time \cite{Be75,In04}. Note that the number of edges incident to some nodes in the resulting graph may be larger than $k$, as illustrated in Figure \ref{2D graph visualization segmentation} where $k=2$, due to symmetry of the undirected graph. The construction of the graph itself provides a basic segmentation of the nodes, for instance in Figure \ref{2D graph visualization segmentation}, it can be observed that the graph contains 3 different connected components. This fact has been exploited in basic graph based classification methods \cite{Al92}. 


In several recent works, the classification problem has been formulated as finding an optimal partition $\{V_i\}_{i=1}^n$ of the nodes $V$ in the graph $G$. The most basic objective function can be formulated as
\begin{align}\label{original problem}
\min_{\{V_i\}_{i=1}^n} & \sum_{i=1}^n \sum_{\substack{(x,y) \in E \; : \\ x \in V_i, \: y \notin V_i}} w(x,y), \\
\text{s.t.} \, & \cup_{i=1}^n V_i  \, = \, V \, , \quad
V_k \cap V_l \, = \, \emptyset \, , \, \forall k \neq l \,
, \label{no_vacuum_overlap}
\end{align}
where the constraint \eqref{no_vacuum_overlap} imposes that there should be no vacuum or overlap between the subsets $\{V_i\}_{i=1}^n$. If $n =2$, then \eqref{original problem} is the so-called "graph cut" \cite{DBLP:conf/coco/1972}. The motivation behind the model \eqref{original problem} is to group the vertices into classes in such a way that pairs of vertices with different class memberships are as dissimilar as possible, indicated by a low value of $w$.


\subsection{Size constraints and supervised constraints}

Extra assumptions are necessary to avoid the trivial global minimizer of \eqref{original problem}, where $V_i = V$ for some $i$ and $V_j = \emptyset$ for all $j \neq i$. There are two common ways to incorporate extra assumptions. In semi-supervised classification problems, the class membership of a small set of the data points is given as training data in advance by the constraints
\begin{equation}\label{supervised constraints}
V_i \supseteq T_i, \quad i \in I = \{1,...,n\},
\end{equation}
where $T_i$ is the set of "training" points known to belong to class $i$. For notational convenience, the set of all class indices $\{1,...,n\}$ is denoted by $I$ in the rest of this paper.

In unsupervised classification problems, one usually assumes that the regions should have approximately equal sizes. The simplest way to achieve this is to impose that each class $V_i$ should have a given size $a_i \in \mathbb{N}$:
\begin{align}\label{constraint balancing}
||V_i|| = a_i, \;  \quad i \in I,
\end{align}
where $\sum_{i=1}^n a_i = ||V||$.
We focus on the case that the norm $||V_i||$ is the number of nodes in $V_i$ for simplicity. As an alternative, $||V_i||$ could be the sum of degrees of each node in $V_i$, where the degree of a node is the number of edges incident to that node. If size constraints are introduced, the problem cannot generally be solved exactly due to NP-hardness. This will be discussed in more detail in Section \ref{Binary representation and convex relaxation of multi-classification problems}.


Usually, a more flexible option is preferred of modifying the energy function such that partitions of equal sizes have lower energy. In case of two classes, the energy \eqref{original problem} becomes $\text{cut}(V_1,V_2)$ $= \sum_{x,y} w(x,y)$, where $x \in V_1$ and $y \in V_2$.
Several different ways of normalizing the energy by the class sizes have been proposed, which can be summarized as follows
\vspace{-0.15cm}
 \begin{equation}\label{balanced cut}
         \text{cut}(V_1,V_2) \Big(\frac{1}{|V_1|} + \frac{1}{|V_2|}\Big), \quad \frac{\text{cut}(V_1,V_2)}{\text{min}(|V_1|,|V_2|)}.
\end{equation}
The expression on the left is called the ratio cut in case of the norm $|V| = \sum_{x \in V}$ and the normalized cut in case of $|V|= \sum_{x \in V} \text{degree}(x)$. The expression on the right is called the Cheeger-ratio cut with the norm $|V| = \sum_{x \in V}$.

The energy functions \eqref{balanced cut} are highly non-convex, but ways to simplify the energy landscape have been proposed \cite{bresson:laurent,szlam,HB10,BTCS14} in order to reduce the number of local minima.

\subsection{New flexible constraint and penalty term on class sizes}

In this paper, we aim to provide a broader set of constraints and penalty terms acting on the class sizes that can be handled in a
completely convex manner. They are designed to achieve the same net result as the ratio energies \eqref{balanced cut} of promoting classes of equal sizes, but in a completely convex way. They can also promote any other size
relations between the class sizes. We will consider flexible size constraints of the form
\begin{align}\label{constraint balancing interval}
S_i^\ell \leq ||V_i|| \leq S_i^u, \;  \quad i \in I,
\end{align}
where $S_i^u \in \mathbb{N}$ is an upper bound on the size of class $i$ and $S_i^\ell \in \mathbb{N}$ is a lower bound. Such types of constraints have previously been proposed for image segmentation in \cite{DBLP:conf/iccv/KlodtC11}. In case one only knows an estimate of the expected class sizes, such constraints can be used to enforce the sizes to lie within some interval of those estimates. To be well defined, it is obviously required that $\sum_{i=1}^n S_i^\ell \leq ||V||$ and $\sum_{i=1}^n S_i^u \geq ||V||$. Note that if $S_i^\ell = S_i^u = a_i$, then \eqref{constraint balancing interval} becomes equivalent to \eqref{constraint balancing}.

\begin{figure}
 \begin{center}
  \begin{tabular}[h!]{c}
   \includegraphics[width=0.3\textwidth]{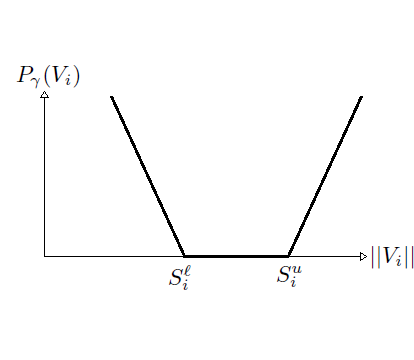}
  \end{tabular}
 \end{center}
 \vspace{-0.5in}
 \caption{\label{penalty term}Illustration of penalty term $P_\gamma(V_i)$.}
\end{figure}

To avoid imposing absolute upper and lower bounds on the class sizes, we also propose appending a piecewise-linear penalty term $\sum_{i=1}^n P_\gamma(||V_i||)$
to the energy function \eqref{original problem}, defined as 
\begin{equation}\label{size penalty term}
P_\gamma(||V_i||) \, = \, \left\{
\begin{array}{ll}
0  \, \quad \; \quad \; \quad \; \quad \text{if} \; S_i^\ell \leq & ||V_i|| \leq S_i^u \\
\gamma \big( ||V_i|| - S_i^u \big) \, \quad \quad \text{if} \; & ||V_i|| > S_i^u \\
\gamma \big( S_i^\ell - ||V_i|| \big) \, \quad \quad \text{if} \; & ||V_i|| < S_i^\ell
\end{array}
\right.
\end{equation}
An illustration of $P_\gamma(||V_i||)$ is given in Figure \ref{penalty term}. In the limit as $\gamma \rightarrow \infty$, the penalty term becomes an equivalent representation of the hard constraints \eqref{constraint balancing interval}. Note that quadratic or higher order penalty terms, although they are convex, are not well suited for the convex relaxation,
because they tend to encourage non-binary values of the labeling functions. In fact, we believe the set
of constraints and penalty terms given here is complete when it comes to being suited for completely convex
relaxations. 

One major contribution of this paper is an efficient algorithmic framework that handles size constraints of the form
\eqref{constraint balancing interval} and the penalty term \eqref{size penalty term} naturally, with almost no additional computational efforts.


\subsection{Region homogeneity terms}

The classification problem \eqref{original problem} involves the minimization of an energy on the boundary of the
classes. The energy is minimized if the data points on each side of the boundary are as dissimilar as possible. These classification models are therefore similar to edge-based image segmentation models, which align the boundary of the regions along edges in the image where the intensity changes sharply. By contrast, region-based image segmentation models, such as the "Chan-Vese" model, use region homogeneity terms that measure how well each pixel fits with each region, in the energy function.

\begin{figure}
 \begin{center}
  \begin{tabular}[t]{c}
   \includegraphics[width=0.45\textwidth]{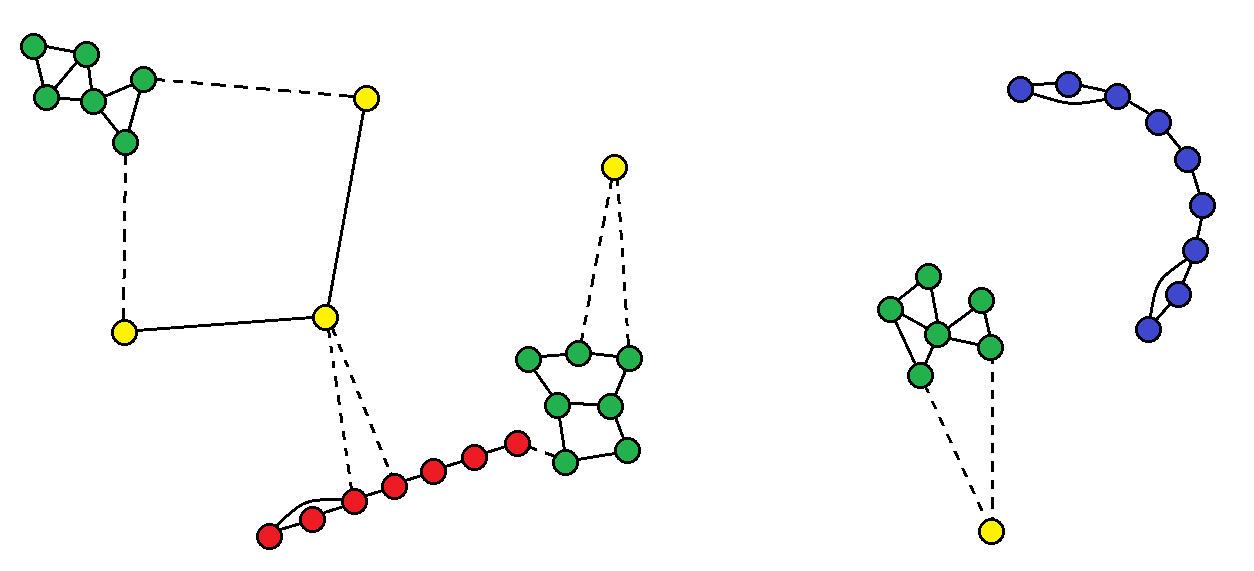}
  \end{tabular}
 \end{center}
 \vspace{-0.2in}
 \caption{\label{2D graph visualization segmentation}Example of segmentation of a graph of 2D points (with number of neighbors $k=2$) into regions of low density (yellow), high degree of correlation of coordinates between neighboring points (red), medium correlation (blue) and low correlation (green). Dashed edges indicate those that contribute to the energy.}
\end{figure}
%

A graph extension of variational segmentation problems with two classes was formulated in \cite{lezoray,LEL12,LEL14,TLE15},
using a non-local total variation term together with a region
term promoting homogeneity of a vertex function. The vertex function could be constructed directly from point geometry and/or from external information such as a color vector defined at each point. We extend the general problem to multiple classes and optional constraints as follows: 



%

\begin{align}\label{original problem with fidelity}
\min_{\{V_i\}_{i=1}^n} \, &  \, \sum_{i=1}^n \sum_{x \in V_i} f_i(x) + \sum_{i=1}^n \sum_{\substack{(x,y) \in E \; : \\ x \in V_i, \: y \notin V_i}} w(x,y) \,,\\
\text{s.t.} \, & \cup_{i=1}^n V_i  \, = \, V \, , \quad
V_k \cap V_l \, = \, \emptyset \, , \, \forall k \neq l \, \nonumber
\end{align}
under optional supervised constraints \eqref{supervised constraints} and/or size constraints
\eqref{constraint balancing interval}/penalty term \eqref{size penalty term}. In \cite{lezoray,LEL12,LEL14,TLE15}
the region terms $f_i(x)$ were defined in terms of a general vertex function $f^0$, which could depend on a color vector or the point geometry. Experimental results on point clouds were shown in case $f^0$ was a
color vector defined at each point. In this work, we will give concrete constructions of $f_i$ for point cloud segmentation purely based on the geometry of the 3D points themselves. For example, the eigenvalues and eigenvectors obtained from a local PCA around each point carry useful information for describing the homogeneity within each class. Concrete examples are given in Section \ref{Segmentation of point clouds}. An illustrative example is given in Figure \ref{2D graph visualization segmentation}, where each node is a 2D point and the region terms have been constructed to distinguish points with different statistical relations to their neighboring points.

The independent preprint \cite{YTS16}, proposed to use region terms in the energy function for semi-supervised classification and the authors proposed a region term that was inferred from the supervised points by diffusion. In contrast, the region terms in this work do not rely on any supervised points, but are as mentioned only specified and demonstrated for the application of 3D point cloud segmentation.

\section{Convex relaxation of minimization problem and analysis based on duality}\label{Binary representation and convex relaxation of multi-classification problems}

In this section, the classification problems are formulated as optimization problems in terms of binary functions instead of sets. The binary representations are used to derive convex relaxations. First, some essential mathematical concepts are introduced, such as various differential operators on graphs. These concepts are used extensively to formulate the binary and convex problems and the algorithms.

\subsection{Differential operators on graphs}\label{Differential operators on graphs}

\indent \indent Our definitions of operators on graphs are based on the theory in \cite{HAL05,ELB08,vangennip:bertozzi}. More information is found in these papers.

Consider two Hilbert spaces, $\mathcal{V}$ and $\mathcal{E}$, which are associated with the sets of vertices and edges, respectively, and the following inner products and norms:
\begin{multline*}
 \hspace{1.1cm}  \langle u,\gamma \rangle_{V} = \sum_{x \in V}u(x) \gamma(x) d(x)^r,  \\   \langle \psi,\phi \rangle_{\mathcal{E}} = \frac{1}{2} \sum_{x,y \in V} \psi(x,y) \phi(x,y) w(x,y)^{2q-1},
\end{multline*}
\begin{equation*}
\norm{u}_{\mathcal{V}} = \sqrt{\langle u, u \rangle_{\mathcal{V}}}= \sqrt{ \sum_{x \in V} u(x)^2d(x)^r},
\end{equation*}
\begin{equation*}
\norm{\phi}_{\mathcal{E}} =   \sqrt{\langle \phi, \phi \rangle_{\mathcal{E}}}= \sqrt{ \frac{1}{2} \sum_{x,y \in V} \phi(x,y)^2 w(x,y)^{2q-1}},
\end{equation*}
\begin{equation}\label{infinity_norm}
\norm{\phi}_{\mathcal{E},\infty} = \max_{x,y \in V}{|\phi(x,y)|},
\end{equation}
for some $r \in [0,1]$ and $q \in [\frac{1}{2},1]$. Above $d(x)$ is the degree of node $x$ (it's number of incident nodes) and
$w(.,.)$ is the weighting function.

From these definitions, we can define the gradient operator $\nabla$ and the Dirichlet energy as
\begin{equation}
    (\nabla u)_w(x,y)= w(x,y)^{1-q}(u(y)-u(x)),
\label{gradient}
\end{equation}
\begin{equation}
    \frac{1}{2} \norm{\nabla u}^2_{\mathcal{E}} = \frac{1}{4} \sum_{x,y \in V} w(x,y) (u(y)-u(x))^2.
\end{equation}

We use the equation $\langle \nabla u, \phi \rangle_{\mathcal{E}} = - \langle u, \diver_w \phi \rangle_{\mathcal{V}}$ to define the divergence:
\begin{equation}
      (\diver_w \phi)(x)= \frac{1}{2d(x)^r} \sum_{y \in V} w(x,y)^q(\phi(x,y)- \phi(y,x)),
\label{divergence}
\end{equation}
where we have exploited symmetry $w(x,y)=w(y,x)$ of the undirected graph in the derivation of the operator.

Using divergence, a family of total variations $TV_{w}:\mathcal{V} \rightarrow \mathbb{R}$ can now be defined:
\begin{multline}\label{tv definition}
     TV_{w}(u) = \sup \big\{ \langle \diver_w \phi, u \rangle_{\mathcal{V}} : \phi \in \mathcal{E}, \norm{\phi}_{\mathcal{E},\infty} \leq 1 \big\} \\
           = \frac{1}{2} \sum_{x,y \in V} w(x,y)^q |u(y)-u(x)|.
\end{multline}

The definition of a family of graph Laplacians $\triangle_w = \diver_w \dot \nabla: \mathcal{V} \rightarrow \mathcal{V}$ is:
\begin{equation}
    (\triangle_w u)(x)= \sum_{y \in V} \frac{w(x,y)}{d(x)^r}(u(y)-u(x)).
\label{Laplacian}
\end{equation}

\subsection{Binary formulation of energy minimization problem}\label{Binary formulation of classification problem}

A partition $\{V_i\}_{i=1}^n$ of $V$ satisfying the no vacuum and overlap constraint
\begin{equation}\label{vac overlap} \cup_{i=1}^n V_i  \, = \, V \, , \quad
V_k \cap V_l \, = \, \emptyset , \, \quad \forall k \neq l
\, \end{equation}
can be described by a binary vector function $u = (u_1,...,u_n) \, : \; V \mapsto \{0,1\}^n$ defined as
\begin{equation}\label{binary partition}
u_i(x) \, := \, \left\{
\begin{array}{ll}
1 , \, & x \in V_i \\
0 , \, & x \notin V_i
\end{array}
\right. \, , \quad i=1,\ldots,n \, .
\end{equation}
In other words, $u(x) = e_i$ if and only if $x \in V_i$, where $e_i$ is the unit normal vector which is $1$ at the $i^{th}$ component and $0$ for all other components. The no vacuum and overlap constraint \eqref{vac overlap} can be expressed in terms of $u$ as
\begin{equation}\label{vac overlap binary}
\sum_{i=1}^n u_i(x) = 1 \, , \quad \forall x \in V.
\end{equation}

Moreover, note that the minimization term of \eqref{original problem} can be naturally related to total variation
\eqref{tv definition} for $q=1$. In fact,
\begin{equation}
\sum_{i=1}^n \sum_{\substack{(x,y) \in E \; : \\ x \in V_i, \: y \notin V_i}} w(x,y) = \sum_{i=1}^n TV_w(u_i).
\end{equation}
This connection between the two terms was used in several recent works to derive, utilizing the graphical framework, efficient
unsupervised algorithms for clustering. For example, \cite{bresson:laurent,TSBLB14} formulate rigorous convergence results for two
methods that solve the relaxed Cheeger cut problem, using non-local total variation. Moreover, \cite{szlam} provides a
continuous relaxation of the Cheeger cut problem, and derives an efficient algorithm for finding good cuts. The authors of \cite{szlam}
relate the Cheeger cut to total variation, and then present a split-Bregman approach of solving the problem. In \cite{TS16} the
continuum limit of total variation on point clouds was derived.

The general set of problems \eqref{original problem with fidelity} can now be formulated in terms of $u$ as
\begin{align}\label{general problem}
\min_{u \in \mathcal{B}}  E^P(u) = \sum_{i=1}^n \sum_{x \in V} C_i(x) u_i(x) + \sum_{i=1}^n TV_w(u_i) 
\end{align}
where
\begin{equation}\label{binary constraint}
\mathcal{B} = \{u : V \mapsto \{0,1\}^n, \quad \sum_{i=1}^n u_i(x) = 1, \;\; \forall x \in V \}
\end{equation}
is the set of binary functions indicating the partition. The superscript $P$ stands for "primal". The optional size
constraints \eqref{constraint balancing interval}, can be imposed in the terms of $u$ as
$$S_i^\ell \leq ||u_i|| \leq S_i^u, \;  \quad i \in I,$$
where $||u_i|| = \sum_{x \in V} u_i(x)$. The size penalty term \eqref{size penalty term} can be imposed by appending the energy function \eqref{general problem} with the term $\sum_{i=1}^n P_{\gamma}(||u_i||)$.

In case of semi-supervised classification, $C_i(x)$ takes the form of
\begin{equation}\label{fidelity function}
C_i(x) = \eta(x) \sum_{i=1}^n |e_i(x) - u^0_i(x)|^2,  
\end{equation}
where $u^0_i$ is a binary function taking the value of $1$ in $T_i$ and $0$ elsewhere, and $\eta(x)$ is a function that takes on a
large constant value $\eta$ on supervised points $\cup_{i=1}^n T_i$ and zero elsewhere. If $\eta$ is chosen sufficiently large, it
can be guaranteed that the solution $u$ satisfies the supervised constraints. 
The algorithm to be presented in this work does not require the selection of an appropriate value for the parameter $\eta$, as
the ideal case where $\eta = \infty$ can be handled naturally without introducing numerical instabilities. 

Region homogeneity terms can be imposed by setting $C_i(x) = f_i(x)$. More generally, region homogeneity terms and supervised data points can be combined by setting
\begin{equation}\label{fidelity function}
C_i(x) = \eta(x) \sum_{i=1}^n |e_i(x) - u^0_i(x)|^2 + f_i(x),  
\end{equation}
The total variation term is defined as in \eqref{tv definition} with $q=1$.

If the number of supervised points is very low and there is no additional region term, the global minimizer of \eqref{general problem} may become the trivial solution where for one of the classes, say $k$, $u_k(x) = 1$ everywhere, and for the other classes $u_i(x) = 1$ for supervised points of class $i$ and $0$ elsewhere. The threshold tends to occur around less than $2.5 \%$ of the points. As in our previous work \cite{MBBT15}, this problem can be countered by increasing the number of edges incident to supervised points in comparison to other points. Doing so will increase the cost of the trivial solution without significantly influencing the desired global minimizer. An alternative, proposed in the preprint \cite{YTS16}, is to create region terms in a pre-processing step by diffusing information of the supervised points into their neighbors.

\subsection{Convex relaxation of energy minimization problem}

Due to the binary constraints \eqref{binary constraint}, the problem \eqref{general problem} is non-convex. As in several recent works on variational image segmentation \cite{ZGFN08,Lellmann-et-al-09a,DBLP:journals/siamis/LiNZS10,BYT2011,BCB10_2,DBLP:journals/jmiv/SawatzkyTJB13} and MRF optimization \cite{DBLP:journals/jacm/KleinbergT02,DBLP:conf/cvpr/BoykovVZ98,Komodakis07approximatelabeling,DBLP:conf/cvpr/AnguelovTCKGHN05}, we replace the indicator constraint set \eqref{binary constraint} by the convex unit simplex
\begin{equation}\label{binary constraint relaxed}
\mathcal{B}' = \{u \, : \: V \mapsto [0,1]^n, \quad \sum_{i=1}^n u_i(x) = 1, \;\; \forall x \in V  \}.
\end{equation}
Hence, we are interested in solving the following convex relaxed problem
\begin{align}\label{general problem g relaxed}
\min_{u \in \mathcal{B}'} E^P(u) = \sum_{i=1}^n \sum_{x \in V} C_i(x) u_i(x) + \sum_{i=1}^n TV_w(u_i).
\end{align}
under optional size constraints \eqref{constraint balancing interval} or penalty term \eqref{size penalty term}.
In case $n=2$ and no size constraints, the relaxation is exact, as proven for image segmentation in \cite{strang08,MR2246072} and classification problems on graphs in \cite{lezoray,MBBT15}. In case $n>2$, the problem becomes equivalent to a multiway cut problem, which is known to be NP-hard \cite{DJPSY92}. In case size constraints are imposed, the problem becomes NP-hard even when $n=2$, as it becomes equivalent to a knapsack problem \cite{martello} in the special case of no TV term.

In this paper we are interested in using the convex relaxation \eqref{general problem g relaxed} to solve the original problem approximately. Under certain conditions, the convex relaxation gives an exact global minimizer of the original problem. For instance, it can be straight forwardly shown that
\begin{proposition}\label{primal: exact global minimizer}
Let $u^*$ be a solution of the relaxed problem \eqref{general problem g relaxed}, with optional size constraints \eqref{constraint balancing interval} or penalty term \eqref{size penalty term}. If $u^* \in \mathcal{B}$, then $u^*$ is a global minimizer of the original non-convex problem \eqref{general problem}.
\end{proposition}
\begin{proof}
Let $E^P(u)$ be the energy function defined in \eqref{general problem g relaxed} with or without the size penalty term \eqref{size penalty term}. Since $\mathcal{B} \subset \mathcal{B}'$ it follows that $\min_{u \in \mathcal{B}'} E^P(u) \leq \min_{u \in \mathcal{B}} E^P(u)$. Therefore, if $u^* = \argmin_{u \in \mathcal{B}'} E^P(u)$ and $u^* \in \mathcal{B}$ it follows that $E(u^*) = \min_{u \in \mathcal{B}} E^P(u)$. The size constraints \eqref{constraint balancing interval} can be regarded as a special case by choosing $\gamma = \infty$.
\end{proof}
If the computed solution of \eqref{general problem g relaxed} is not completely binary, one way to obtain an approximate binary solution that exactly indicates the class membership of each point, is to select the binary function as the nearest vertex in
the unit simplex by the threshold
\begin{equation}\label{rounding}u^T(x) = e_\ell(x), \; \text{where} \; \ell
= \argmax_{i \in I} u_i(x).
\end{equation}
As an alternative to the threshold scheme \eqref{rounding}, binary solutions of the convex relaxation
\eqref{general problem g relaxed} can also be obtained from a dual solution of \eqref{general problem g relaxed},
which has a more solid theoretical justification if some conditions are fulfilled. The dual problem also gives insight
into why the convex relaxation is expected to closely approximate the original problem. This is the topic of the next section.

\subsection{Analysis of convex relaxation through a dual formulation}\label{analysis dual}

We will now derive theoretical results which indicate that the multiclass problem \eqref{general problem}
is closely approximated by the convex relaxation \eqref{general problem g relaxed}. The following results extend those given in \cite{BYT2011} from image domains to graphs. In contrast to \cite{BYT2011}, we also incorporate size constraints or penalty terms in the analysis. In fact, the strongest results given near the end of the section are only valid for problems without such size constraints/terms. This observation agrees well with our experiments, although in both cases very close approximations are obtained.

We start by deriving an equivalent dual formulation of \eqref{general problem g relaxed}. Note that this dual problem is different from the "max-flow" type dual problem on graphs proposed
in our previous work \cite{MBBT15} in case of two classes. Its main purpose is theoretical analysis,
not algorithmic development. In fact, its relation to flow maximization will be the subject of the next section.
Dual formulations on graphs have also been proposed in \cite{HLE13} for variational multiscale decomposition of graph signals.

\begin{theorem}\label{prop dual problem}
The convex relaxed problem \eqref{general problem g relaxed} can equivalently be formulated as the dual problem
\begin{align}\label{dual problem}
\sup_{q,\rho^1,\rho^2} \sum_{x \in V} 
\min_{i\in I}& \Big(C_i(x) + (\diver_w q_i)(x) + \rho_i^2 - \rho_i^1 \Big)\nonumber \\
 + & \big( \rho_i^1 S^\ell_i - \rho_i^2 S^u_i \big),
\end{align}
subject to
\begin{align}
& (q_1,...,q_n) \in S^n_\infty, \label{q-constraint} \\ 
& \rho_i^1,\rho_i^2 \in [0,\gamma], \quad i=1,...,n, \label{rho constraints}
\end{align}
where the above set of infinity norm spheres is defined as
\begin{align}\label{infinity norm spheres}
& S^n_\infty = \{(q_1,....,q_n) \, : \mathcal{E} \mapsto \mathbb{R}^n \; \text{s.t.} \; \norm{q_i}_{\mathcal{E},\infty} \leq 1 \, \forall i\}.
\end{align}
No size information is incorporated by choosing $\gamma = 0$. The size penalty term
\eqref{size penalty term} is incorporated by choosing $0 < \gamma < \infty$. Size constraints \eqref{constraint balancing interval}
are incorporated by choosing $\gamma = \infty$.

\end{theorem}
\begin{proof}
By using the definition of total variation \eqref{tv definition}, the problem \eqref{general problem g relaxed} with size penalty term
\eqref{size penalty term} can be expressed in primal-dual form as
\begin{align} \label{eqn:two-opt}
& \min_{u \in \mathcal{B}'} \,  \sup_{q \in S^n_\infty } \; \sum_{i=1}^n P(||u_i||) \nonumber \\
& + \sum_{i=1}^n \sum_{x \in V}  u_i(x) \big(C_i(x) \, +\, (\diver_w q_i)(x) \big),
\end{align}
where $S^n_\infty$ is defined in \eqref{infinity norm spheres}. 
It will be shown that the size constraints \eqref{constraint balancing interval} or penalty term \eqref{size penalty term} can be implicitly incorporated by introducing the dual variables $\rho_i^1,\rho_i^2 \in \mathbb{R}_+$, $i=1,..,n$ as
\begin{align} \label{eqn:two-opt-2}
& \min_{u \in \mathcal{B}'} \,  \sup_{q \in S^n_\infty, \rho^1,\rho^2 \in [0,\gamma]^n } E(u;q,\rho^1,\rho^2) \nonumber \\
& = \sum_{i=1}^n \sum_{x \in V}  u_i(x) \big\{C_i(x) \, +\, (\diver_w q_i)(x) + \rho_i^2 - \rho_i^1 \big\} \nonumber \\
& \quad +  \big( \rho_i^1 S^\ell_i - \rho_i^2 S^u_i \big),
\end{align}
The primal-dual problem \eqref{eqn:two-opt-2} satisfies all the conditions of the mini-max theorem (see e.g. Chapter 6,
Proposition 2.4 of \cite{ET99}). The constraint sets for $q,\rho^1,\rho^2$ and $u$ are compact and convex, and the energy
function $E(u,q)$ is convex l.s.c. for fixed $q$ and concave u.s.c. for fixed $u$. This implies the existence of at
least one primal-dual solution (saddle point) of finite energy value.

For a given $u$, the terms involving $\rho^1$ and $\rho^2$ can be rearranged as
\begin{align}\label{balancing term-n}
& \sup_{ 0 \leq \rho_i^1 \leq \gamma}   \rho^1_i \big( S^\ell_i - \sum_{x \in V} u_i(x) ) & \nonumber \\
& = \left\{
\begin{array}{ll}
0 \,  & \mbox{ if  } \sum_{x \in V}  u_i(x) \geq S^\ell_i \\
\gamma\big(S^\ell_i-\sum_{x \in V}u_i(x)\big) \,  &  \text{ if  }  \sum_{x \in V}  u_i(x) < S^\ell_i
\end{array}\right.
\end{align}
\begin{align}\label{balancing term 2-n}
& \sup_{0 \leq \rho_i^2 \leq \gamma }  \rho_i^2 \big( \sum_{x \in V} u_i(x) - S^u_i ) & \nonumber \\
& = \left\{
\begin{array}{ll}
0 \,  & \mbox{ if  } \sum_{x \in V}  u_i(x) \leq  S^u_i \\
\gamma \big(\sum_{x \in V} u_i(x)-S^u_i\big) \,  &  \text{ if  }  \sum_{x \in V}  u_i(x) >  S^u_i
\end{array}\right. &
\end{align}
Consider the above three choices for $\gamma$. In case $\gamma = 0$ the class sizes do not contribute to the energy. In case $0 < \gamma < \infty$
the two above terms summed together is exactly equal to the size penalty term $P(||u_i||)$. In case $\gamma = \infty$, the constraint
set on $\rho^1,\rho^2$ is no longer compact, but we can apply Sion's generalization of the mini-max theorem \cite{Si58},
which allows either the primal or dual constraint set to be non-compact. It follows that if the size constraints \eqref{constraint balancing interval} are not satisfied,
the energy would be infinite, contradicting existence of a primal-dual solution.

From the mini-max theorems, it also follows
that the inf and sup operators can be interchanged as follows
\begin{align}
& \min_{u \in \mathcal{B}'} \, \sup_{q \in S^n_\infty, \rho^1,\rho^2 \in [0,\gamma]^n} \; E(u;q,\rho^1,\rho^2) \, = \, \nonumber \\
& \sup_{q \in S^n_\infty, \rho^1,\rho^2 \in [0,\gamma]^n} \,\min_{u \in \mathcal{B}'} \; E(u;q,\rho^1,\rho^2). \label{eqn:minimax}
\end{align}
For notational convenience, we denote the unit simplex pointwise as
\begin{equation}\label{unit simplex pointwise}
\Delta^n_+ = \{(u_1,...,u_n) \in [0,1]^n \, : \; \sum_{i=1}^n u_i = 1 \}
\end{equation}
For an arbitrary vector $F = (F_1, \ldots, F_n) \in
\mathbb{R}^n$, observe that  
\begin{equation} \label{eq:u-indicator}
\min_{(u_1, \ldots, u_n) \in \triangle_{+}} \, \sum_{i=1}^n u_i F_i
\, = \, \min(F_1, \ldots,F_n) \,.
\end{equation}
Therefore, the inner minimization of \eqref{eqn:minimax} can be solved analytically at each position $x \in V$,
and we obtain the dual problem
\begin{align}
& \sup_{q \in S^n_\infty, \rho^1,\rho^2 \in [0,\gamma]^n} \big( \rho_i^1 S^\ell_i - \rho_i^2 S^u_i \big) \nonumber \\
& + \sum_{x \in V} \min_{u(x) \in \triangle_{+}} \sum_{i=1}^n u_i \big\{ C_i \, +\, \diver_w q_i + \rho_i^2 - \rho_i^1 \big\}(x) \nonumber \\
& = \sup_{q \in S^n_\infty, \rho^1,\rho^2 \in [0,\gamma]^n} \big( \rho_i^1 S^\ell_i - \rho_i^2 S^u_i \big) \nonumber \\
& + \sum_{x \in V} \min_{i \in I}\{C_i(x) + (\diver_w q_i)(x) + \rho_i^2 - \rho_i^1 \}. \nonumber
\end{align}

\end{proof}

Assuming a solution of the dual problem $q^*,{\rho^1}^*,$ ${\rho^2}^*$ has been obtained, the following theorem characterizes
the corresponding primal variable $u^*$
\begin{theorem} \label{theo:consist_unq}
There exists a maximizer  $q^{\ast},{\rho^1}^*,{\rho^2}^*$ to the dual problem \eqref{dual problem}. At the point $x \in V$, let $I_m(x) = \{i_1,...,i_k\}$
be the set of indices such that
\begin{equation}\label{minimal component}
I_m(x) = \argmin_{i \in I}\Big( C_i(x) + (\diver_w q^*_i)(x) + {\rho_i^2}^* - {\rho_i^1}^* \Big).
\end{equation}
There exists a solution $u^*$ to the primal problem \eqref{general problem g relaxed} such that $(u^*;q^*,{\rho^1}^*,{\rho^2}^*)$ is a primal-dual pair. At the point $x$, $u^*(x)$ must satisfy
\begin{equation} \label{eq:cons02}
\sum_{i \in I_m(x)} u_{i}^{\ast}(x) \, = \, 1 \; \text{ and} \quad u_j^{\ast}(x) \,
= \, 0 \, , \, j \notin I_{\min} \, .
\end{equation}
If the minimizer \eqref{minimal component} is unique at the point $x \in V$, then the corresponding primal solution $u^*$ at the point $x$ must be valued
\begin{equation}\label{binary primal solution}
u^*_i(x) = \, \left\{
\begin{array}{ll}
1 , \, & \text{if}\;i = I_m(x) \\
0 , \, & \text{if}\;i \neq I_m(x)
\end{array}
\right. \, , \quad i=1,\ldots,n \, .
\end{equation}
If the minimizer \eqref{minimal component} is unique at every point $x \in V$, then the corresponding primal solution $u^*$, given by the formula \eqref{binary primal solution}, is an exact global binary minimizer of the original non-convex problem \eqref{general problem}.
\end{theorem}
\begin{proof}
Since all conditions of the mini-max theorem \cite{ET99,Si58} are satisfied (c.f. proof of Theorem \ref{prop dual problem}), there must exist a maximizer $q^*,{\rho^1}^*,{\rho^2}^*$ of the dual problem \eqref{dual problem} and a minimizer $u^*$ of the primal problem \eqref{general problem g relaxed} such that
$(u^*,q^*)$ is a solution of the primal-dual problem \eqref{eqn:two-opt} (see e.g. \cite{ET99}). For arbitrary vectors
$u \in \Delta_+^n$ and $F \in \mathbb{R}^n$, it must hold that $\sum_{i\in I} u_i F_i \geq \min_{i \in I} F_i$. Therefore, at the point $x$, $u^*$ must satisfy
\begin{align}
& \sum_{i\in I} u^*_i(x)\big( (C_i + \diver_w q^*_i)(x) + {\rho_i^2}^* - {\rho_i^1}^* \big) \nonumber \\
& = \min_{i \in I} \big( (C_i + \diver_w q_i)(x) + {\rho_i^2}^* - {\rho_i^1}^* \big), \nonumber 
\end{align}
otherwise the primal-dual energy would exceed the dual energy, contradicting strong duality. The above expression can be further decomposed as follows
\begin{align*}
& = \sum_{i\in I_m(x)} u^*_i(x)\big(  (C_i + \diver_w q^*_i)(x)  + {\rho_i^2}^* - {\rho_i^1}^* \big)\\
& + \sum_{i\notin I_m(x)} u^*_i(x)\big(  (C_i + \diver_w q^*_i)(x) + {\rho_i^2}^* - {\rho_i^1}^* \big) \\
& = \big( \sum_{i\in I_m(x)} u^*_i(x)\big)\min_{i \in I}\big( (C_i + \diver_w q^*_i)(x)  + {\rho_i^2}^* - {\rho_i^1}^* \big)\\
& + \sum_{i\notin I_m(x)} u^*_i(x)\big(  (C_i + \diver_w q^*_i)(x)  + {\rho_i^2}^* - {\rho_i^1}^* \big)
\end{align*}
Since $\big(  (C_j + \diver_w q^*_j)(x)  + {\rho_i^2}^* - {\rho_i^1}^*\big)(x)$\\
 $> \min_{i \in I}\big( (C_i + \diver_w q^*_i)(x)  + {\rho_i^2}^* - {\rho_i^1}^*(x)$ for all $j \notin I_m$, the above can only be true provided
$\sum_{i \in I_m} u^*_i = 1$ and $u_i^*(x) = 0$ for $i \notin I_m$.

If the minimizer $I_m(x)$ is unique, it follows directly from \eqref{eq:cons02}, that $u^*_i(x)$ must be the indicator vector \eqref{binary primal solution}.

If the minimizer $I_m(x)$ is unique at every point $x \in V$, then the corresponding primal solution $u^*$ given by \eqref{binary primal solution} is contained in the binary set $\mathcal{B}$. By Proposition \ref{primal: exact global minimizer}, $u^*$ is a global minimizer of \eqref{general problem}.
\end{proof}

%
%

It can also be shown that an exact binary primal solution exists if there are two non-unique minimal components to the vector
\begin{equation*}
(C(x) + \diver_w q^{\ast}(x) + {\rho^2}^* - {\rho^1}^*) 
\end{equation*}
but this result only holds in case there are no constraints acting on the class sizes.
\begin{theorem}\label{theorem two components}
Assume that $q^*$ is a maximizer of the dual problem \eqref{dual problem} with $\gamma = 0$, i.e. no class size constraints. If
\eqref{minimal component} has at most two minimal components for all $x \in V$,
then there exists a corresponding binary primal solution to the convex relaxed primal problem \eqref{general problem g relaxed},
which is a global minimizer of the original non-convex problem \eqref{general problem}.
\end{theorem}
A constructive proof of Theorem \ref{theorem two components} is given in Appendix \ref{appendix 1}.

If the vector $(C(x) + \diver_w q^{\ast}(x) + {\rho^2}^* - {\rho^1}^*)$ has three or more minimal components, it cannot in general be expected that a corresponding binary primal solution exists, reflecting that one can probably not obtain an exact solution to the NP-hard problem \eqref{general problem} in general by a convex relaxation. Experiments indicate that this very rarely, if ever, happens in practice for the classification problem \eqref{general problem}.

As an alternative thresholding scheme, $u^T$ can be selected based on the formula \eqref{binary primal solution} after a dual solution to the convex relaxation has been obtained. If there are multiple minimal components to the vector $(C + \diver q^*)(x)$, one can select $u^T(x)$ to be one for an arbitrary one of those indices, just as for the ordinary thresholding scheme \eqref{rounding}. 
Experiments will demonstrate and compare both schemes in Section \ref{Applications and experiments}.

\section{\lq Max-flow' formulation of dual problem and algorithm}\label{Max-flow dual formulations and algorithms}

A drawback of the dual model \eqref{dual problem} is the non-smoothness of the objective function, which is also a drawback of the original primal formulation of the convex relaxation. This section reformulates the dual model in a structurally similar way to a max-flow problem, which is smooth and facilitates the development of a very efficient algorithm based on the augmented Lagrangian theory.

The resulting dual problem can be seen as a multiclass variant of the max-flow model proposed in our work \cite{MBBT15} for two classes, and a graph analogue of the max-flow model given for image domains in \cite{YBTB10}. Note that our derivations differ conceptually from \cite{YBTB10,MBBT15}, because we directly utilize the dual problem derived in the last section. Furthermore, the new flexible size constraint \eqref{constraint balancing interval} and penalty term \eqref{size penalty term} are incorporated naturally in the max-flow problem
by a modified flow conservation condition, which indicates that there should be a constant flow
excess at each node. The amount of flow excess is expressed with a few additional optimization variables in the algorithm,
and they can optimized over with very little additional computational cost.

\subsection{\lq Max-flow' reformulation dual problem}\label{Without balancing constraints}

We now derive alternative dual and primal-dual formulations of the convex relaxed problem that are more beneficial for computations. The algorithm will be presented in the next section.
\begin{proposition}\label{dual to max-flow}
The dual problem \eqref{dual problem} can equivalently be formulated as the dual \lq max-flow' problem:
\begin{equation}
\label{eq:max-flowc2} \sup_{p_s, p, q,\rho^1,\rho^2} \;
\,\sum_{x \in V} \, p_s(x) + \sum_{i=1}^n \big( \rho_i^1 S^\ell_i - \rho_i^2 S^u_i \big)\, \end{equation} 
subject to, for all $i \in I$,
\begin{align}
&  \abs{q_i(x,y)} \,  \leq \, 1,  \quad \quad \quad \quad \quad \quad \quad \;\;\, \forall  (x,y)  \in  E, \label{eq:cond-01c2}\\
& p_i(x) \,  \leq  \, C_i(x), \quad \quad \quad \quad \quad \quad \quad \quad \;\;\;\;\, \forall x \in  V,  \label{eq:cond-03c2}\\
& \big( \diver_w q_i - p_s + p_i \big)(x) \, = \, \rho_i^1 - \rho_i^2, \quad \; \forall x \in  V, \label{eq:cond-04c2} \\
& 0 \leq \rho_i^1,\rho_i^2 \leq \gamma. \label{eq:cond-05c2}
\end{align}
\end{proposition}
\begin{proof}
By introducing the auxiliary variable $p_s \, : \; V \mapsto \mathbb{R}$, the dual problem \eqref{dual problem} can be reformulated as follows
\begin{align}
&\sup_{q, \rho^1,\rho^2 } \, \sum_{x \in V} \min_{i\in I} \big( C_i(x) + \diver_w q_i(x)  + \rho_i^2 - \rho_i^1 ) \big) \,\nonumber \\
& \quad \quad + \sum_{i=1}^n \big( \rho_i^1 S^\ell_i - \rho_i^2 S^u_i \big) \nonumber \\
& \text{subject to, for all $i \in I$,}\nonumber \\
& \norm{q_i}_{\mathcal{E},\infty} \leq 1, \nonumber \\
& 0 \leq \rho_i^1,\rho_i^2 \leq \gamma. \nonumber
\end{align}
\begin{align}
& = \sup_{p_s, q, \rho^1,\rho^2}
\sum_{x \in V} p_s(x) + \sum_{i=1}^n \big( \rho_i^1 S^\ell_i - \rho_i^2 S^u_i \big)\, \nonumber \\ 
& \text{subject to, for all $i \in I$,}\nonumber \\
& p_s(x) \leq  (C_i + \diver_w q_i)(x) + \rho_i^2 - \rho_i^1 \;\; \forall x \in V, \label{constraints} \\
& \norm{q_i}_{\mathcal{E},\infty} \leq 1, \nonumber \\
& 0 \leq \rho_i^1,\rho_i^2 \leq \gamma.\nonumber 
\end{align}
By adding another set of auxiliary variables $p_i \, : \; V \mapsto \mathbb{R}$, $i =1,...,n$, the constraints \eqref{constraints} can be formulated as
\begin{align}
& p_s(x) = p_i(x) + \diver_w q_i(x) + \rho_i^2 - \rho_i^1,  \label{flow conservation}\\
& p_i(x) \leq C_i(x), \nonumber
\end{align}
for all $x \in V$ and all $i \in I$. Rearranging the terms in constraint \eqref{flow conservation}, and using the definition of the infinity norm \eqref{infinity_norm}, leads to the \lq max-flow' model \eqref{eq:max-flowc2} subject to \eqref{eq:cond-01c2}-\eqref{eq:cond-05c2}.
\end{proof}
Problem \eqref{eq:max-flowc2} with constraints \eqref{eq:cond-01c2}-\eqref{eq:cond-05c2} is structurally similar to a max-flow problem over $n$ copies of the graph $G$, $(V_1,E_1) \times ... \times (V_n,E_n)$, where $(V_i,E_i) = G$ for $i \in I$. The aim of the max-flow problem is to maximize
the flow from a source vertex to a sink vertex under flow capacity at each edge and flow conservation at each node. The variable
$p_s(x)$ can be regarded as the flow on the edges from the source to the vertex $x$ in each of the subgraphs $(V_1,E_1), ..., (V_n,E_n)$, which have unbounded capacities. The variables $p_i(x)$ and $C_i(x)$ can be regarded as the flow and capacity on
the edge from vertex $x$ in the subgraph $(V_i,E_i)$ to the sink. Constraint \eqref{flow conservation} is the flow conservation condition. Observe that in case of size constraints/terms, instead of being conserved, there should be a constant excess flow $\rho_i^1 - \rho_i^2$ for each node in the subgraph $(V_i,E_i)$. The objective function \eqref{eq:max-flowc2} is a measure of the total amount of flow in the graph.

Utilizing results from Section \ref{analysis dual}, we now show that the convex relaxation \eqref{general problem g relaxed} is the equivalent dual problem to the max-flow problem \eqref{eq:max-flowc2}.
\begin{theorem}\label{dual-primaldual-primal}
The following problems are equivalent to each other:\\
1) The max-flow problem \eqref{eq:max-flowc2}, subject to \eqref{eq:cond-01c2}-\eqref{eq:cond-05c2};\\\\
2) The primal-dual problem:
\begin{align} \label{eq:primal-dual}
& \min_{u}  \sup_{p_s, p, q,\rho^1,\rho^2}  \Big\{E(p_s,
p, q, \rho^1,\rho^2; u) \nonumber \\
 & = \sum_{x \in V} p_s(x) + \sum_{i=1}^n \big( \rho_i^1 S^\ell_i - \rho_i^2 S^u_i \big) \nonumber \\
 & + \sum_{i=1}^n \sum_{x \in V} u_i \big(\diver_w q_i - p_s + p_i + \rho_i^2 - \rho_i^1 \big)(x) \Big\}
 \end{align}
subject to \eqref{eq:cond-01c2}, \eqref{eq:cond-03c2} and \eqref{eq:cond-05c2}, where $u$ is the relaxed region indicator function.\\

3) The convex relaxed problem \eqref{general problem g relaxed} with size constraint \eqref{constraint balancing interval} if $\gamma = \infty$, size penalty term \eqref{size penalty term} if $0 < \gamma < \infty$ and
no size constraints if $\gamma = 0$.
\end{theorem}
\begin{proof}
The equivalence between the primal-dual problem \eqref{eq:primal-dual} and the max-flow problem \eqref{eq:max-flowc2} follows directly as $u_i$ is an unconstrained Lagrange multiplier for the flow conservation constraint \eqref{flow conservation}. Existence of the Lagrange multipliers follows as: 1) \eqref{eq:max-flowc2} is upper bounded, since it is equivalent to \eqref{dual problem}, which by Theorem \ref{theo:consist_unq}
admits a solution; 2) the constraints \eqref{eq:cond-04c2} are linear, and hence differentiable.

The equivalence between the primal-dual problem \eqref{eq:primal-dual}, the max-flow problem \eqref{eq:max-flowc2} and the convex relaxed problem \eqref{general problem g relaxed} now follows: By Proposition \ref{dual to max-flow} the \lq max-flow' problem \eqref{eq:max-flowc2} is equivalent to the dual problem \eqref{dual problem}. By Theorem \ref{prop dual problem}, the dual problem \eqref{dual problem} is equivalent to the convex relaxed problem \eqref{general problem g relaxed}
with size constraints \eqref{constraint balancing interval} if $\gamma = \infty$, size penalty term \eqref{size penalty term} if $0 < \gamma < \infty$ and no size constraints if $\gamma = 0$.

\end{proof}

Note an important distinction between the primal-dual problem \eqref{eq:primal-dual} and the primal-dual problem \eqref{eqn:minimax} derived in the last section: The primal variable $u$ is unconstrained in \eqref{eq:primal-dual}. The simplex constraint $\mathcal{B}'$ is handled implicitly. It may not seem obvious from the proof how the constraints on $u$ are encoded in the primal-dual problem, therefore we give some further insights: For a given primal
variable $u$, the maximization with respect to $p_s$ of the primal-dual problem \eqref{eq:primal-dual} at the point $x$ can be rearranged as
\begin{align}\label{opt p_s} \sup_{p_s(x)} & ((1- \sum_{i=1}^n  u_i )p_s)(x) \nonumber \\ = & \left\{
\begin{array}{ll}
0 \,  & \mbox{ if  } \sum_{i=1}^n  u_i(x) = 1 \\
\infty \,  &  \text{ if  }  \sum_{i=1}^n  u_i(x) \neq 1
\end{array}\right.
\end{align}
If $u$ does not satisfy the sum to one constraint at $x$, then the primal-dual energy would be infinite, contradicting boundedness from above. In a similar manner, the optimization with respect to $p_i$ can be expressed as
\begin{align} \label{opt
p_t} \sup_{p_i(x) \leq C_i(x)} & u_i(x) p_i(x) = \left \{
\begin{array}{ll}
(u_i C_i)(x) \,  & \; \mbox{ if  } u_i(x) \geq 0 \\ \infty \,
& \; \mbox{ if  } u_i(x) < 0 .
\end{array}\right.
\end{align}
which would be infinite if $u$ does not satisfy the non-negativity constraints. If $u(x) = e_i$, the indicator function of class $i$, the value would be $C_i(x)$, which is indeed the pointwise cost of assigning $x$ to class $i$.

\subsection{Augmented Lagrangian max-flow algorithm}\label{sec algorithm}

This section derives an efficient algorithm, which
exploits the fact all constraints on $u$ are handled implicitly in the primal-dual problem \eqref{eq:primal-dual}.
The algorithm is based on the augmented Lagrangian theory, where $u$ is updated as a Lagrange multiplier by a gradient descent
step each iteration. Since no subsequent projection of $u$ is necessary, the algorithm tolerates a wide range of step sizes and converges with high accuracy. The advantages of
related \lq max-flow' algorithms for ordinary 2D imaging problems over e.g. Arrow-Hurwicz type primal-dual
algorithms have been demonstrated in \cite{DBLP:conf/emmcvpr/BaeTY14,YBTB2013}.

From the primal-dual problem \eqref{eq:primal-dual}, we first construct the augmented Lagrangian functional:

\begin{align}
& L_{c} = \sum_{x \in V} p_{s} +  \sum_{i=1}^n ( \rho_i^1 S^\ell_i - \rho_i^2 S^u_i ) \big) \nonumber \\
&\hspace{0.5cm}+ \sum_{x \in V} u_i(x) \big( \diver_w q_i  - p_{s}+ p_{i} + \rho_i^2 - \rho_i^1 \big)(x) \nonumber
\\&\hspace{0.5cm} - \frac{c}{2} \sum_{i=1}^n \norm{\diver_w q_i - p_{s}+ p_{i} + \rho_i^2 - \rho_i^1}_2^2. \label{augmented lagrangian}
\end{align}
An augmented Lagrangian algorithm for minimizing the above functional is given below, which involves alternatively maximizing $L_c$ for the dual variables and then updating the Lagrange multiplier $u$.

Note that if there are no constraints on the class sizes, $\gamma = 0$, then ${\rho^1}^k = {\rho^2}^k \equiv 0$ for every iteration $k$. The algorithm can in this case be simplified by setting ${\rho^1}^k = {\rho^2}^k \equiv 0$ for all $k$ and ignoring all steps involving $\rho^1$ and $\rho^2$.

\line(1,0){218}

\fbox{\textbf{Algorithm 1} }\\ 


Initialize $p_s^1$, $p^1$, $q^1$, ${\rho^1}^1$, ${\rho^2}^1$ and $u^1$. For $k=1,...$ until convergence:
\begin{itemize}
\renewcommand\labelitemi{$\bullet$}
\item Optimize $q$ flow, for $i \in I$
\begin{multline}
\label{updatep_s}
\hspace{-0.5cm}
q_i^{k+1}  =  \argmax_{|q(e)| \leq
1 \; \forall e \in E}
-\frac{c}{2}
 \norm{\diver_w q - F^k}_2^2,
\end{multline}
where  $F^k= {p_s}^{k}  - {p_i}^{k} + \frac{u_i^k}{c} - {\rho_i^2}^k + {\rho_i^1}^k$ is fixed.

\item Optimize source flow $p_s$
\begin{multline}
\label{updateps2}
\hspace{-0.5cm}
p_s^{k+1} =   \argmax_{p_s(x) }  \sum_{x \in V} \,
 p_s  - \frac{c}{2}\norm{p_s - G^k}_2^2,
\end{multline}
where $G^k= {p_i}^k + \diver_w q_i^{k+1} -  \frac{u^k_i}{c} + {\rho_i^2}^k - {\rho_i^1}^{k+1}$ is fixed.

\item Optimize sink flow $p_i$, for $i \in I$,
\begin{multline}
\label{updatept}
p_i^{k+1}  := \argmax_{p_i(x) \leq \; C_i(x) \; \forall x \in V} \,
-  \frac{c}{2}\norm{p_i - H^k }_2^2,
\end{multline}
where $H^k= {p_s}^{k+1} - \diver_w {q}_i^{k+1} + \frac{u^k_i}{c} - {\rho_i^2}^k + {\rho_i^1}^k$ is fixed.

\item Optimize $\rho^1_i$, for $i \in I$,
\begin{multline}
\label{rho1}
{\rho^1_i}^{k+1}  =  \argmax_{0 \leq \rho^1_i \leq \gamma} \sum_{x \in V} \rho_i^1 S^\ell_i -\frac{c}{2} \norm{\rho_i^1 - J^k}_2^{2},
\end{multline}
where  $J^k =  - p^{k+1}_{i} -  \mbox{div}_w q_i^{k+1} + \frac{u^{k}_i}{c} +  p^{k+1}_s - {\rho^2}^{k}_i$ is fixed.
\item Optimize $\rho^2_i$, for $i \in I$,
\begin{multline}
\label{rho2}
{\rho^2_i}^{k+1}  =  \argmax_{0 \leq \rho^2_i \leq \gamma} \; \sum_{x \in V} - \rho_i^2 S^u_i -\frac{c}{2} \norm{{\rho_i^2} - M^k}_2^{2},
\end{multline}
where  $M^k =  p^{k+1}_{i} +  \mbox{div}_w q_i^{k+1} - \frac{u^{k}_i}{c} -  p^{k+1}_s - {\rho^1_i}^{k+1}$ is fixed.

\item Update $u_i$, for $i \in I$
\begin{align*}
& u^{k+1}_i \, = \, u^k_i \\
& - c \, (\diver_w q_i^{k+1} - p_s^{k+1} +
p_i^{k+1}+ {\rho^2_i}^{k+1} - {\rho^1_i}^{k+1}).
\end{align*}
\end{itemize}
\line(1,0){218}


The optimization problem \eqref{updatep_s} can be solved by a few steps of the projected gradient method as follows:
\begin{equation}\label{q iteration}
     q^{k+1}_i= \Pi_{W} (q_i + c \nabla_w(\mbox{div}_w q^{k}_i - F^k)),
\end{equation}
Above, $\Pi_{w}$ is a projection operator which is defined as
         \begin{multline}\label{projection operator}
\Pi_{W} (s(x,y))= \\
  \begin{cases}  s(x,y) &\mbox{if } |s(x,y)| \leq 1, \\
      \text{sgn}(s(x,y)) & \mbox{if } |s(x,y)| > 1,
        \end{cases}
\end{multline}
where sgn is the sign function. There are extended convergence theories
for the augmented Lagrangian method in the case when one of the subproblems is solved inexactly, see e.g. \cite{Esser10,GBO09}. In our experience, one gradient ascent iteration leads to the fastest overall speed of convergence.


The subproblems \eqref{updateps2} and \eqref{updatept} can be solved by
\begin{equation}
    p_s(x)= G^k(x)+ \frac{1}{c},
\end{equation}
\begin{equation}
    p_i(x)= \min (H^k(x), C_i(x)).
\end{equation}
Consider now the subproblems (\ref{rho1}) and (\ref{rho2}). In case no constraints are given on $\rho^1$ and $\rho^2$, the maximizers over the sum of the concave quadratic terms can be computed as the average of the maximizers to each individual term as
$$
\mbox{mean} \big( - J^k + \frac{S^\ell_i }{c \,||V||}\big), \quad \mbox{mean} \big( - M^k - \frac{S_i^u}{c\,||V||}  \big),
$$
respectively for $\rho^1$ and $\rho^2$. Since the objective function is concave, and the maximization variable is just a constant, an exact solution to the constrained maximization problem can now be obtained by a projection onto that constraint as follows
\begin{align}
   & {\rho_i^1}^{k+1}= \min\Big(\max\big( \mbox{mean}( - J^k + \frac{S^\ell_i }{c\,||V||}), 0 \big),\gamma\Big), \label{update rho size penalty}\\
   & {\rho_i^2}^{k+1}= \min\Big(\max\big(  \mbox{mean}( - M^k - \frac{S_i^u}{c\,||V||}), 0 \big),\gamma\Big). \label{update rho size penalty 2}
   \end{align}
Algorithm 1 is suitable for parallel implementation on GPU, since the subproblems at each substep can be solved pointwise independently
of each other using simple floating point arithmetics. The update formula \eqref{q iteration} for $q$, which only requires access to the values of neighboring nodes at the previous iterate. As discussed in Section \ref{Segmentation as optimization problem over graph}, the number of neighbors may vary for nodes across the graph, therefore special considerations should be taken when declaring memory. We have implemented the algorithm on CPU for experimental evaluation for simplicity.

\section{Applications and experiments}\label{Applications and experiments}

We now focus on specific applications of the convex framework. Experimental results on semi-supervised classification of high-dimensional data are presented in Section \ref{Semi-supervised Classification Results}. Section \ref{Segmentation of point clouds} proposes specific terms in the general model \eqref{original problem with fidelity} for segmentation of unstructured 3D point clouds, and presents experimental results on LaDAR data acquired in outdoor scenes.
In both cases we give a thorough examination of accuracy of the results, tightness of the convex relaxations, and convergence properties of the algorithms.

A useful quality measure of the convex relaxation is to what extent the computed solution is binary. Proposition \ref{primal: exact global minimizer} indicates that if the computed solution is completely binary, it is also an exact global minimizer to the original non-convex problem. Let ${u^k}^T$ be a thresholding of $u^k(x)$ in the sense that each row of $u^k$ is modified to be the closest vertex in the unit simplex according to the scheme \eqref{rounding}. As a quality measure of the solution $u^k$ at each iteration $k$ of Algorithm 1, we calculate the average difference between $u^k$ and its thresholded version ${u^k}^T$ as follows:
%
%
%
\begin{equation}\label{binary_difference_formula}
      b(u^k) = \frac{1}{2n||V||} \big(\sum_{i=1}^n \sum_{x \in V} |{u^k}^T_i(x) - u^k_i(x))|\big),
\end{equation}
where $||V||$ is the number of nodes, and $n$ is the number of classes. We call $b(u^k)$ the "binary difference" of $u^k$ at
iteration $k$. Naturally, we want $b(u^k)$ to become as low as possible
as the algorithm converges.

\subsection{Semi-supervised classification results}\label{Semi-supervised Classification Results}

\indent \indent In this section, we describe the supervised classification results, using the algorithm with and without the size constraints \eqref{constraint balancing}, \eqref{constraint balancing interval} or penalty term \eqref{size penalty term}.

We compare the accuracy of the results with respect to the ground truth. The results are also compared against other local minimization
approaches in terms of the final total variation energies:
\begin{equation*}
      E(u) = \frac{1}{2} \sum_{i=1}^n \sum_{x,y \in V} w(x,y)|u_i(x) - u_i(y))|,
\end{equation*}
where $n$ is the number of classes. A lower value of $E$ is better. The energy contribution from the fidelity term is ignored because the solution satisfies the supervised constraints by construction, thus giving zero contribution from those terms.

To compute the weights for the data sets, we use the Zelnik-Manor and Perona weight function \cite{perona}. The function is defined as:
\begin{equation}
        w(x,y)=\exp\left(-\frac{d(x,y)^{2}}{\sqrt{\tau(x)\tau(y)}}\right),
\label{perona_eq}
\end{equation}
where $d(x,y)$ is a distance measure between vertices $x$ and $y$, and $\sqrt{\tau(x)}$ is the distance between vertex $x$ and its $M^{th}$ closest neighbor.  If $y$ is not among the $M$ nearest neighbors of $x$, then $w(x,y)$ is set to 0. After the graph is computed, we symmetrize it by setting
\begin{equation*}
w(x,y)=$ $\max(w(x,y),w(y,x)).
\end{equation*}
Here, $M$ is a parameter to be chosen. The weight function will be defined more specifically for each application.

We run the minimization procedure until the following stopping criterion is satisfied:
\begin{equation*}
     \frac{1}{||V||} \big(\sum_i \sum_{x \in V} |u_i(x) - u^{old}_i(x))|\big) < \delta,
\end{equation*}
where $u^{old}$ is the $u$ from the previous iteration, and the value of $\delta$ varies depending on the data set (anywhere from $10^{-12}$ to $10^{-10}$).

All experiments were performed on a 2.4 GHz Intel Core i2 Quad CPU. We initialize $C_i(x)=$ constant (in our case, the constant is set to 500) if $x$ is a supervised point of any class but class $i$, and $0$ otherwise, for all $i \in I$. The variables $u$, $q_i$, $\rho_i^1$, $\rho_i^2$ are initialized to zero for all $i \in I$. The variable $p_s$ is initialized to $C_n$, where $n$ is the number of classes. We set $p_i=p_s \hspace{0.15cm} \forall i \in I$.


In the following, we give details about the set up and results for each dataset, before we draw some general conclusions in the end.

\subsubsection{MNIST}
\begin{figure}[h]
\centering
  \scalebox{1.5}  {\includegraphics[trim=50 40 50 20,clip,width=2in,height=1in]{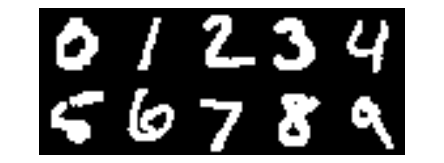}}
\caption{Examples of digits from the MNIST data base}
\label{fig:digits}

\end{figure}

\indent \indent The MNIST data set~\cite{lecun:cortes}, affiliated with the Courant Institute of New York University, consists of $70,000$ $28 \times 28$ images of handwritten digits $0$ through $9$. Some of the images in the database are shown in Figure \ref{fig:digits}. The objective is, of course, to assign the correct digit to each image; thus, this is a $10$-class segmentation problem.

We construct the graph as follows; each image is a node on a graph, described by the feature vector of $784$ pixel intensity values in the image. These feature vectors are used to compute the weights for pairs of nodes. The weight matrix is computed using the Zelnik-Manor and Perona weight function \eqref{perona_eq} with local scaling using the $8^{th}$ closest neighbor. We note that preprocessing of the data is not needed to obtain an accurate classification; we do not perform any preprocessing. The parameter $c$ used was 0.05.

The average accuracy results over 100 different runs with randomly chosen supervised points
are shown in Table \ref{table:MNIST} in case of no size constraints.
We note that the new approaches reach consistently higher accuracies and lower energies than related local minimization approaches, and that incorporation of size information can improve the accuracies further. The computation times are highly efficient, but not quite as fast as MBO, which only uses 10 iteration to solve the problem in an approximate manner. The $\text{Log}_{10}$ plots of the binary difference versus iteration, depicted in Figure \ref{fig:binary_difference}, show that the binary difference converges to an extremely small number.

\begin{figure*}[t]
  \subfigure[ground truth]{\scalebox{0.65}{\includegraphics[trim=50 30 10 30,clip,width=3.25in,height=2in]{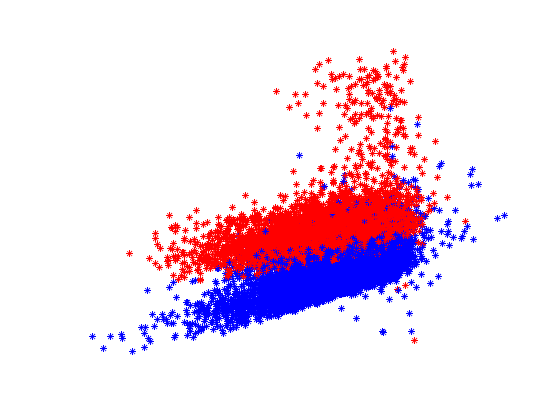}}}
        \subfigure[Proposed result (randomly selected supervised points)]{ \scalebox{0.65}{       \includegraphics[trim=50 30 10 30,clip,width=3.25in,height=2in]{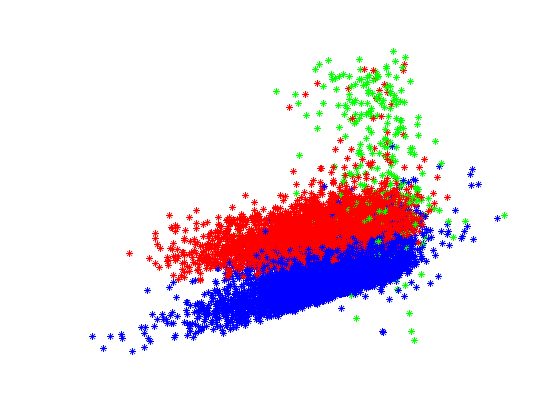}}}
         \subfigure[Proposed result (non-randomly selected supervised points)]{\scalebox{0.65}{\includegraphics[trim=50 30 10 30,clip,width=3.25in,height=2in]{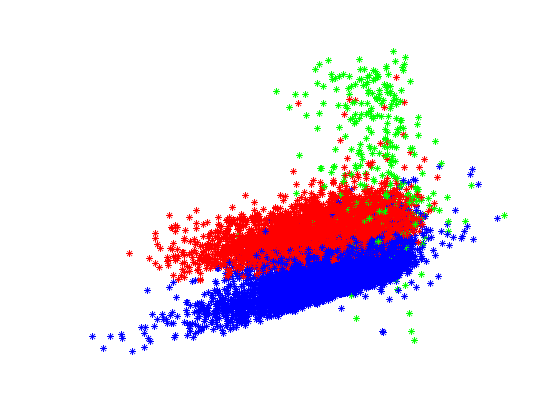}}}
\begin{center}
        \subfigure[MBO result (randomly selected supervised points)]{ \scalebox{0.65}{       \includegraphics[trim=50 30 10 30,clip,width=3.25in,height=2in]{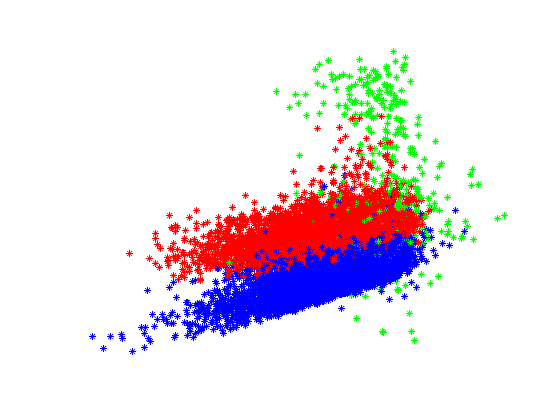}}}
     \subfigure[MBO result (non-randomly selected supervised points)]{ \scalebox{0.65} {\includegraphics[trim=50 30 10 30,clip,width=3.25in,height=2in]{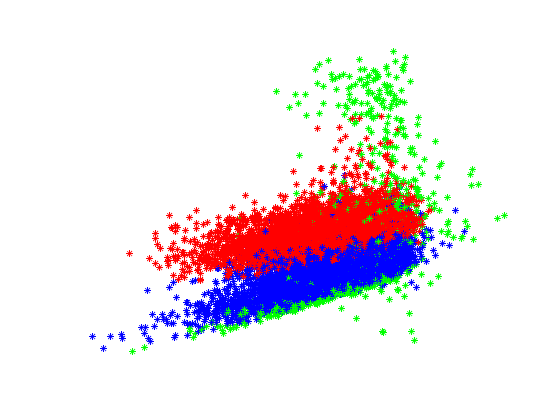}}}
\end{center}
\caption{MNIST results. These graphs visualize the values of the first versus the sixth eigenvector (of the graph Laplacian) relating to the nodes of class 4 and 9 only. The blue and red region represents nodes of class 4 and 9, respectively. The green region represents misclassified points.}
\label{fig:mnist}
\end{figure*}

The results of the data set are visualized in Figure \ref{fig:mnist}. For the visualization procedure, we use the first and the sixth eigenvector of the graph Laplacian. The dimension of each of the eigenvectors is $N \times 1$, and each node of the data set is associated with a value of each of the vectors. One way to visualize a classification of a data set such as MNIST, which consists of a collection of images, is to plot the values of one eigenvector of the graph Laplacian versus another and use colors to differentiate classes in a given segmentation. In this case, the plots in Figure \ref{fig:mnist} graph the values of the first versus the sixth eigenvector (of the graph Laplacian) relating to the nodes of class 4 and 9 only. The blue and red region represents nodes of class 4 and 9, respectively. The green region represents misclassified points.

Moreover, we compare our results to those of other methods in Table \ref{benchmark_main}, where our method's name is written in bold. Note that algorithms such as linear and nonlinear classifiers, boosted stumps, support vector machines and both neural and convolution nets are all supervised learning approaches, which use around $60,000$ of the images as a training set ($86\%$ of the data) and $10,000$ images for testing. However, we use only $3.57\%$ (or less) of our data as supervised training points, and obtain classification results that are either competitive or better than those of some of the best methods. Moreover, note that no preprocessing was performed on the data, as was needed for some of the methods we compare with; we worked with the raw data directly.

\subsubsection{Three Moons Data Set}

\indent \indent We created a synthetic data set, called the three moons data set, to test our method.
The set is constructed as follows. First, consider three half circles in $\mathbb{R}^2$.
The first two half top circles are unit circles with centers
at $(0,0)$ and $(3,0)$. The third half circle is a bottom half circle with radius of $1.5$ and
center at $(1.5,0.4)$. A thousand points from each of the three half
circles are sampled and embedded in $\mathbb{R}^{100}$ by adding
Gaussian noise with standard deviation of $0.14$ to each of the $100$
components of each embedded point. The goal is to segment the circles, using
a small number of supervised points from each class. Thus, this is a 3-class segmentation problem. The noise and the fact that the points are embedded in high-dimensional space make this difficult.

We construct the graph as follows; each point is a node on a graph, described by the feature vector consisting of the $100$ dimensions of the point. To compute the distance component of the weight function for pairs of nodes, we use these feature vectors. The weight matrix is computed using the Zelnik-Manor and Perona weight function \eqref{perona_eq} with local scaling using the $10^{th}$ nearest neighbor. The parameter $c$ was 0.1.

\begin{figure}[t]
         \subfigure[ground truth]{\scalebox{0.45}{\includegraphics[trim=50 30 10 30,clip,width=3.25in,height=2in]{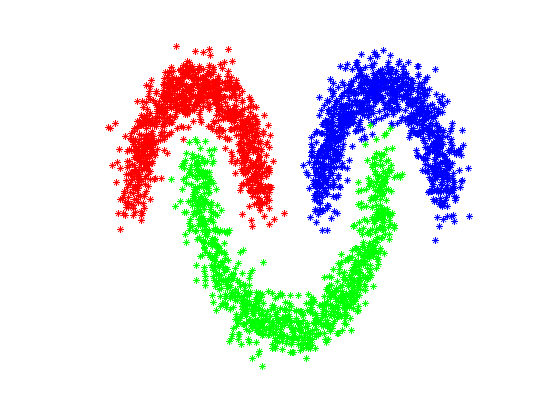}}}
        \subfigure[Proposed result]{ \scalebox{0.45}{       \includegraphics[trim=50 30 10 30,clip,width=3.25in,height=2in]{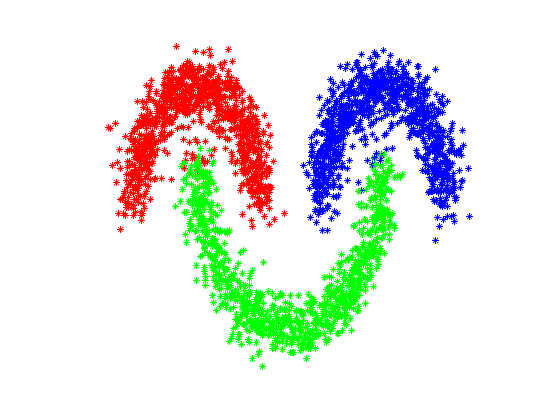}}}
\begin{center}
     \subfigure[MBO result]{ \scalebox{0.45} {\includegraphics[trim=50 30 10 30,clip,width=3.25in,height=2in]
{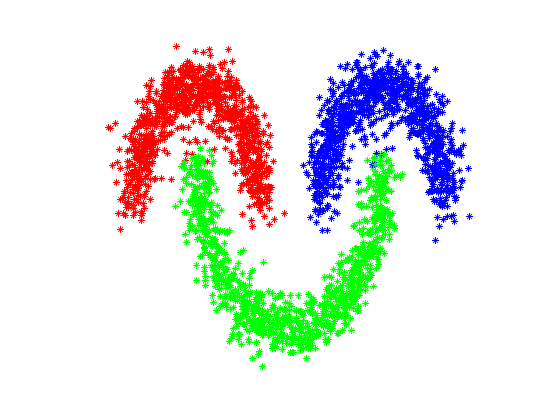}}}
\end{center}
\caption{Three moons results}
\label{fig:moons}
\end{figure}

The results of the data set are visualized in Figure \ref{fig:moons} and the accuracies are shown in Table \ref{benchmark_main}. This is the only dataset where the proposed approach got lower accuracy than MBO. For this particular example, the global minimizer does not seem the best in terms of accuracy, which is a fault of the model rather than an optimization procedure.

\begin{figure*}[t]
     \subfigure[ground truth]{\scalebox{0.65}{\includegraphics[trim=50 30 10 30,clip,width=3.25in,height=2in]{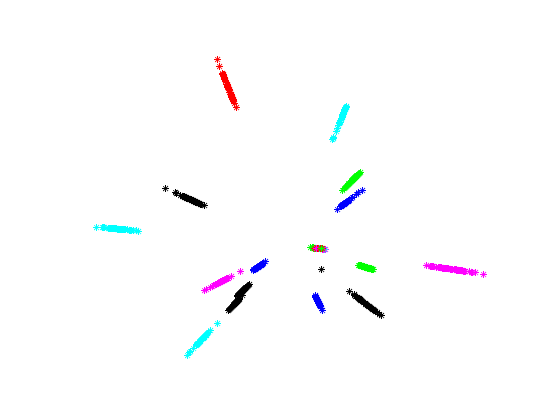}}}
 \subfigure[MBO result (non-randomly selected supervised points)]{ \scalebox{0.65} {\includegraphics[trim=50 30 10 30,clip,width=3.25in,height=2in]{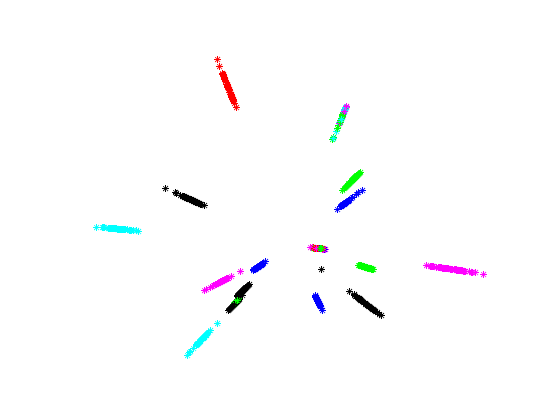}}}
         \subfigure[Proposed result (non-randomly selected supervised points)]{\scalebox{0.65}{\includegraphics[trim=50 30 10 30,clip,width=3.25in,height=2in]{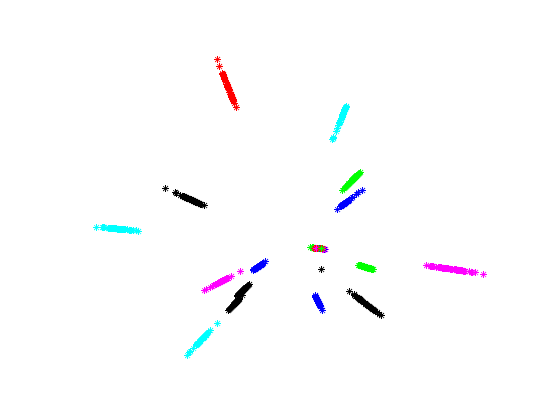}}}
\caption{COIL Results. These graphs visualize the values of the first versus the third eigenvector of the graph Laplacian.
The results of the classification are labeled by different colors.}
\label{fig:coil}
\end{figure*}


   \subsubsection{COIL}
     \indent \indent We evaluated our performance on the benchmark COIL data set~\cite{coil,chapelle:scholkopf:zien} from the Columbia University Image Library. This is a set of color $128 \times 128$  images of $100$ objects, taken at different angles. The red channel of each image was downsampled to $16 \times 16$ pixels by averaging over blocks of $8 \times 8$ pixels. Then, $24$ of the objects were randomly selected and then partitioned into six classes. Discarding $38$ images from each class leaves $250$ per class, giving a data set of $1500$ data points and 6 classes.

We construct the graph as follows; each image is a node on a graph. We apply PCA to project each image onto 241 principal components; these components form the feature vectors. The vectors are used to calculate the distance component of the weight function. The weight matrix is computed using the Zelnik-Manor and Perona weight function \eqref{perona_eq} with local scaling using the $4^{th}$ nearest neighbor. The parameter $c$ used was 0.03.

Resulting accuracies are shown in Table \ref{benchmark_main}, indicating that our method outperforms local minimization approaches and is comparable to or better than some of the other best existing methods. The results of the data set are visualized in Figure \ref{fig:coil}; the procedure used is similar to that of the MNIST data set visualization procedure. The plots in the figure graph the values of the first versus the third eigenvector of the graph Laplacian. The results of the classification are labeled by different colors.


\subsubsection{Landsat Satellite data set}

 We also evaluated our performance on the Landsat Satellite data set, obtained from the UCI Machine Learning Repository \cite{uci}. This is a hyperspectral data set which is composed of sets of multi-spectral values of pixels in 3 $\times$ 3 neighborhoods in a satellite image; the portions of the electromagnetic spectrum covered include near-infrared. The goal is to predict the classification of the central pixel in each element of the data set. The six classes are red soil, cotton crop, grey soil, damp grey soil, soil with vegetation stubble and very damp grey soil. There are 6435 nodes in the data set.

   We construct the graph as follows. The UCI website provides a $36$-dimensional feature vector for each node. The feature vectors are used to calculate the distance component of the weight function. The weight matrix is computed using the Zelnik-Manor and Perona weight function \eqref{perona_eq} with local scaling using the $4^{th}$ nearest neighbor. The parameter c used was 0.3.


Table \ref{benchmark_main} includes comparison of our method to some of the best methods (most cited in \cite{mroueh}). One can see that our results are of higher accuracy. We now note that, except the GL and MBO algorithms, all other algorithms we compare the Landsat satellite data to are supervised learning methods, which use 80\% of data for training and 20\% for testing. Our method was able to outperform these algorithms while using a very small percentage of the data set (10\%) as supervised points. Even with 5.6\% supervised points it outperforms all but one of the aforementioned methods.


\subsubsection{Non-uniform distribution of supervised points}

   \indent \indent
   In all previous experiments, the supervised points have been sampled randomly from all the datapoints. To test the algorithms in more challenging scenarios, we introduce some bias in the sampling of the supervised points, which is also a more realistic situation in practice. We used two different data sets for this test: the MNIST data set and the COIL data set.

    In the case of the MNIST data set, we chose the supervised points non-randomly for digits $4$ and $9$ only. To obtain the non-randomness, we allowed a point to be chosen as supervised only if it had a particular range of values for the second eigenvector. This resulted in a biased distribution of the supervised points. The results for this experiment were the following: for the max flow algorithm, the overall accuracy was $97.734\%$, while for digits $4$ and $9$, it was $96.85\%$. For comparison, the non-convex MBO algorithm \cite{garcia} gave an accuracy of $95.60\%$ overall, but $89.71\%$ for digits $4$ and $9$. The MBO method was also a bit more unstable in its accuracy with respect to different distributions of the supervised points. The max-flow algorithm was very stable, with a very small standard deviation for a set of accuracies for different supervised point distributions.

    In the case of the COIL data set, we chose the supervised points non-randomly for classes $2$ and $6$. The non-randomness was achieved in the same way as for the MNIST data set. The results were the following: the overall accuracy of the max-flow was $92.69\%$, while for classes $2$ and $6$, it was $90.89\%$. The MBO algorithm \cite{garcia} gave an accuracy of $83.90\%$ overall, but $77.24\%$ for classes $2$ and $6$.

    These results are summarized in Table \ref{table:non_random} and are visualized in Figures \ref{fig:mnist} and \ref{fig:coil} for MNIST and COIL data sets, respectively.

\begin{figure}[t!]
        \subfigure[algorithm without size constraints]{ \scalebox{1.1}{       \includegraphics[width=2.75in,height=2in]{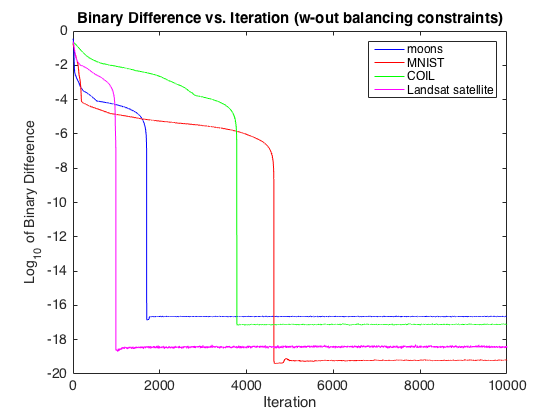}}}
        \subfigure[algorithm with flexible constraints \eqref{constraint balancing interval} and penalty term \eqref{size penalty term} acting on class sizes]{ \scalebox{1.1}{       \includegraphics[width=2.75in,height=2in]{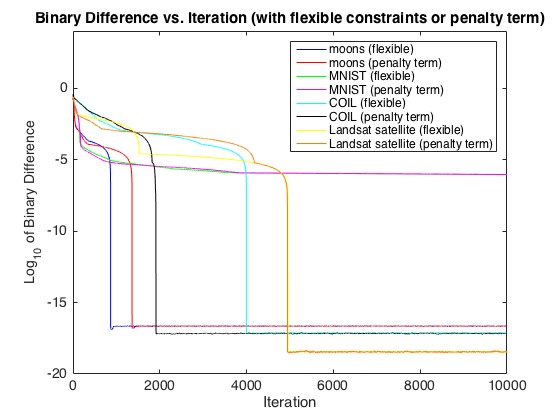}}}
\caption{$\text{Log}_{10}$ plot of binary difference $b(u^k)$ vs. iteration count.}
\label{fig:binary_difference}
\end{figure}

\subsubsection{Experiments with size constraints and penalty term}\label{Experiments with size constraints and penalty term}


The exact size constraints \eqref{constraint balancing} could improve the accuracies if knowledge of the exact class sizes are available.
However, it is not realistic to obtain the exact knowledge of the class sizes in practice, and this was the motivation behind developing the flexible constraints \eqref{constraint balancing interval} or the penalty term \eqref{size penalty term}. In order to simulate the case that only an estimate of the class sizes are known, we perturb the exact class sizes by a random number ranging between 1 $\%$ and 20 $\%$ of $||V||/n$.  The lower and upper bounds in \eqref{constraint balancing interval} and \eqref{size penalty term} are centered around the perturbed class size, and the difference between them is chosen based on the uncertainty of the estimation, which we assume to be known. More specifically, denoting the exact class size $c_i$, the perturbed class size $\tilde{c}_i$ is chosen as a random number in the interval $[c_i-p,c_i+p]$. In experiments, we select $p$ as 1 $\%$, 10 $\%$ and 20 $\%$ of $||V||/n$. 
The lower and upper bounds in the flexible size constraint \eqref{constraint balancing interval} and the penalty term \eqref{size penalty term} are chosen as $S^\ell_i = \tilde{c}_i - p$ and $S^u_i = \tilde{c}_i +p$. The parameter $\gamma$ in the penalty term is set to $10$ for all datasets.

We run the algorithm for each choice of $p$ several times with different random selections of the perturbed class size $\tilde{c}_i$ each time. The average accuracies over all the runs for each choice of $p$ are shown in Table \ref{perturbed}. The flexible size constraints or penalty term improve the accuracy compared to the case when no size information was given, shown in Table \ref{benchmark_main}. Note that the accuracy improves also in cases of great uncertainty in the class size estimates ($p = 20 \%$). The exact size constraints can be seen to not be suitable in case knowledge of the exact class sizes are not known, as imposing them significantly reduces the accuracies in those cases.

\begin{table}[p!]
\begin{minipage}{0.48\textwidth}
\caption{Accuracy compared to ground truth of the proposed algorithm vs. other algorithms.}
\label{benchmark_main}
\vspace{-0.2in}
\begin{center}
\begin{subtable}{}
\caption*{$\quad \quad \quad \quad \quad \quad $\fontsize{9}{9}\selectfont {MNIST (10 classes)\\
Note that some of the comparable algorithms, marked by *, use substantially more data for training ($85.7 \%$ at most and $21.4 \%$ at smallest) than the proposed algorithm, see the main text for more information.} \label{table:MNIST}}
\begin{tabular}{|c| c|}
\hline
Method & Accuracy \\
\hline
p-Laplacian~\cite{buhler} & 87.1\% \\
\hline
multicut normalized 1-cut~\cite{hein} & 87.64\%\\
\hline
linear classifiers~\cite{lecun,lecun:cortes} & 88\% \\
\hline
Cheeger cuts~\cite{szlam} & 88.2\% \\
\hline
boosted stumps* ~\cite{kegl,lecun:cortes} & 92.3-98.74\% \\
\hline
transductive classification~\cite{szlam:maggioni:coifman} & 92.6\% \\
\hline
tree GL~\cite{gc:flenner:percus} & 93.0\% \\
\hline
$k$-nearest neighbors* ~\cite{lecun,lecun:cortes} & 95.0-97.17\% \\
\hline
neural/conv. nets* ~\cite{lecun,ciresan,lecun:cortes} & 95.3-99.65\% \\
\hline
nonlinear classifiers* ~\cite{lecun,lecun:cortes} & 96.4-96.7\% \\
\hline
SVM* ~\cite{lecun,decoste} & 98.6-99.32\% \\
\hline
GL \cite{garcia} (3.57\% supervised pts.) & 96.8\% \\
\hline
MBO \cite{garcia} (3.57\% supervised pts.) & 96.91\% \\
\hline
{\bf Proposed} (3.57\% supervised pts.) & \bf{97.709\%} \\
\hline
\end{tabular}
\end{subtable}
\end{center}

%

\begin{center}
\begin{subtable}{}
\caption*{\fontsize{9}{9}\selectfont {Three moons (5\% supervised points)}\label{table:moons}}
\begin{tabular}{|c| c|}
\hline
Method & Accuracy \\
\hline
GL \cite{garcia} & 98.4\% \\
\hline
MBO \cite{garcia}  & 99.12\% \\
\hline
{\bf Proposed}  & \bf{98.714\%} \\
\hline
\end{tabular}
\end{subtable}
\end{center}

\begin{center}
\begin{subtable}{}
\caption*{\fontsize{9}{9}\selectfont {COIL ($10 \%$ supervised points)} \label{table:COIL}}
\begin{tabular}{|c| c|}
\hline
Method & Accuracy \\
\hline
$k$-nearest neighbors~\cite{subramanya} & 83.5\% \\
\hline
LapRLS~\cite{belkin,subramanya} & 87.8\% \\
\hline
sGT~\cite{joachims,subramanya} & 89.9\% \\
\hline
SQ-Loss-I~\cite{subramanya} & 90.9\% \\
\hline
MP~\cite{subramanya} & 91.1\% \\
\hline
GL \cite{garcia} & 91.2\% \\
\hline
MBO \cite{garcia} & 91.46\% \\
\hline
{\bf Proposed}  & \bf{93.302\%} \\
\hline
\end{tabular}
\end{subtable}

\begin{subtable}{}
\caption*{\fontsize{9}{9}\selectfont {$\quad \quad \quad \quad \quad$ Landsat satellite data set\\
*-marked use $80 \%$ of the data set for training, see main text for more information}\label{table:satellite}}
\begin{tabular}{|c| c|}
\hline
Method & Accuracy \\
SC-SVM* \cite{mroueh} & 65.15\% \\
\hline
SH- SVM* \cite{mroueh} & 75.43\% \\
\hline
S-LS* \cite{mroueh} )& 65.88\% \\
\hline
simplex boosting* \cite{mroueh}  & 86.65\% \\
\hline
S-LS rbf.* \cite{mroueh} & 90.15\% \\
\hline
GL  \cite{garcia} (10\% supervised pts.) & 87.62\%  \\
\hline
GL \cite{garcia} (5.6\% supervised pts.) & 87.05\% \\
\hline
MBO \cite{garcia} (10\% supervised pts.) & 87.76\%  \\
\hline
MBO \cite{garcia} (5.6\% supervised pts.)  & 87.25\% \\
\hline
{\bf Proposed} (10\% supervised pts.) & \bf{90.267\%} \\
\hline
{\bf Proposed} (5.6\% supervised pts.) &  \bf{88.621\%} \\
\hline
\end{tabular}
\end{subtable}
\end{center}

\end{minipage}\hfill
\end{table}
\begin{table}[p!]
\begin{minipage}{0.48\textwidth}
%
\begin{center}
\caption{Accuracies in case of non-uniformly distributed supervised points}
\label{table:non_random}
\begin{tabular}{|c|c|c|}
\hline
 & overall &  classes 4 and 9 \\
\hline
\textbf{Proposed}, MNIST & \bf{97.734\%} & \bf{96.85\%} \\
\hline
MBO, MNIST & 95.60\% & 89.72\% \\
\hline
\hline
 & overall &  classes 2 and 6 \\
\hline
\textbf{Proposed}, COIL & \bf{92.69\%} & \bf{90.89\%} \\
\hline
MBO, COIL & 83.90\% & 77.24\% \\
\hline
\end{tabular}
\end{center}

\caption{Accuracies for experiments with class size incorporation. The exact class sizes are perturbed by a random number within $p$ $\%$ of the size and the accuracies are computed by averaging over multiple runs. See Section \ref{Experiments with size constraints and penalty term} for details.}
\label{perturbed}
\begin{center}
\begin{subtable}{MNIST, 3.57\% supervised points}
\begin{tabular}{|c|c|c|c|}
\hline
max size perturbation ($p$)& $1 \%$ & $10 \%$& $20 \%$ \\
\hline
flexible size constraints \eqref{constraint balancing interval}& 97.761 & 97.725 & 97.716 \\
\hline
penalty term \eqref{size penalty term} & 97.755 & 97.739 & 97.722 \\
\hline
exact size constraints \eqref{constraint balancing} & 96.139 & 70.820 & 63.660 \\
\hline
\end{tabular}
\end{subtable}
\end{center}

\vspace{-0.7cm}

\begin{center}
\begin{subtable}{}
\caption*{\fontsize{9}{9}\selectfont {Three moons}\label{table:moons}}
\begin{tabular}{|c|c|c|c|}
\hline
\multicolumn{4}{|c|}{ $5\%$ supervised points }\\
\hline
max size perturbation ($p$)& $1 \%$ & $10 \%$& $20 \%$ \\
\hline
flexible size constraints \eqref{constraint balancing interval} & 99.374 & 98.829 & 98.750 \\
\hline
penalty term \eqref{size penalty term} & 99.368 & 98.789 & 98.718 \\
\hline
exact size constraints \eqref{constraint balancing} & 99.108 & 72.685 & 66.627 \\
\hline
\multicolumn{4}{|c|}{ $0.6\%$ supervised points }\\
\hline
max size perturbation ($p$)& $1 \%$ & $10 \%$& $20 \%$ \\
\hline
flexible size constraints \eqref{constraint balancing interval} & 97.833 & 97.738 & 97.160 \\
\hline
penalty term \eqref{size penalty term} & 97.848 & 97.793 & 97.406 \\
\hline
exact size constraints \eqref{constraint balancing} & 97.706 & 68.956 & 66.872 \\
\hline
\end{tabular}
\end{subtable}
\end{center}

\vspace{-0.7cm}

\begin{center}
\begin{subtable}{}
\caption*{\fontsize{9}{9}\selectfont {COIL}\label{table:moons}}
\begin{tabular}{|c|c|c|c|}
\hline
\multicolumn{4}{|c|}{ $10\%$ supervised points }\\
\hline
max size perturbation ($p$)& $1 \%$ & $10 \%$& $20 \%$ \\
\hline
flexible size constraints \eqref{constraint balancing interval} & 93.403 & 93.535 & 93.527 \\
\hline
penalty term \eqref{size penalty term} & 93.360 & 93.418 & 93.325 \\
\hline
exact size constraints \eqref{constraint balancing} & 92.990 & 59.936 & 55.624 \\
\hline
\multicolumn{4}{|c|}{ $5\%$ supervised points }\\
\hline
max size perturbation ($p$)& $1 \%$ & $10 \%$& $20 \%$ \\
\hline
flexible size constraints \eqref{constraint balancing interval} & 90.428 & 90.892 & 90.730 \\
\hline
penalty term \eqref{size penalty term} & 89.957 & 90.967 & 90.712 \\
\hline
exact size constraints \eqref{constraint balancing} & 89.931 & 55.152 & 54.674 \\
\hline
\end{tabular}
\end{subtable}
\end{center}

\vspace{-0.7cm}

\begin{center}
\begin{subtable}{}
\caption*{\fontsize{9}{9}\selectfont {Landsat satellite data set }\label{table:moons}}
\begin{tabular}{|c|c|c|c|}
\hline
\multicolumn{4}{|c|}{ $10\%$ supervised points }\\
 \hline
max size perturbation ($p$)& $1 \%$ & $10 \%$& $20 \%$ \\
\hline
flexible size constraints \eqref{constraint balancing interval} & 90.504 & 90.397 & 90.344 \\
\hline
penalty term \eqref{size penalty term} & 90.479 & 90.371 & 90.347 \\
\hline
exact size constraints \eqref{constraint balancing}& 87.773 & 67.687 & 65.757 \\
\hline
\multicolumn{4}{|c|}{ $5\%$ supervised points }\\
  \hline
max size perturbation ($p$)& $1 \%$ & $10 \%$& $20 \%$ \\
\hline
flexible size constraints \eqref{constraint balancing interval} & 89.024 & 89.022 & 88.848 \\
\hline
penalty term \eqref{size penalty term} & 89.025 & 89.018 & 88.987 \\
\hline
exact size constraints \eqref{constraint balancing}& 86.327 & 60.904 & 51.276 \\
\hline
\end{tabular}
\end{subtable}
\end{center}

\end{minipage}

\end{table}


\begin{table}[t]
\caption{Timing results (in seconds)}
\label{table:timing}
\begin{center}
\begin{tabular}{|c|c|c|c|}
\hline
  & MBO \cite{garcia} & GL \cite{garcia} & Proposed \\
\hline
MNIST & 15.4 & 153.1 & \bf{42.5} \\
\hline
3 moons & 3.7 & 3.9  & \bf{2.7} \\
\hline
COIL & 1.18 & 1.19 &  \bf{1.4} \\
\hline
satellite & 16.4 & 23 & \bf{16.5} \\
\hline
\end{tabular}
\end{center}
\end{table}

\begin{table}[t]
\caption{Initial and Final Energy}
\label{table:energy}
\begin{center}
\begin{tabular}{|c|c|c|c|}
\hline
& initial energy & final energy & final energy \\
& & (MBO) \cite{garcia} & proposed \\
\hline
MNIST & 225654 & 15196 & \bf{12324} \\ 
\hline
3 moons & 5982.79 &  433.19 & \bf{420.24}\\
\hline
COIL & 1774.3 & 24.61 & \bf{24.18} \\
\hline
satellite & 5116.9 & 221.87 & \bf{214.95} \\
\hline
\end{tabular}
\end{center}
\end{table}

\subsubsection{Summary of experimental results}

Experimental results on the benchmark datasets, shown in Table \ref{benchmark_main}, indicate a consistently higher accuracy of the proposed convex algorithm than related local minimization approaches based on the MBO or Ginzburg-Landau scheme. The improvements are especially significant when the supervised points are not uniformly distributed among the dataset as shown in Table \ref{table:non_random}. On one synthetic dataset, "three moons", the accuracy of the new algorithm was slightly worse, indicating that the global
minimizer was not the best in terms of accuracy for this particular toy example. Table \ref{table:energy} shows that the new algorithm reaches the lowest
energy on all of the experiments, further indicating that MBO and Ginzburg-Landau are not able to converge to the global minimum.
Table \ref{benchmark_main} shows that the accuracies of the proposed algorithm are also highly competitive against a wide range of other established algorithms, even when substantially less training data than those algorithms are being used. Table \ref{perturbed} shows that that the flexible size constraints \eqref{constraint balancing interval} and penalty term \eqref{size penalty term} can improve the accuracy, if a rough estimate of the approximate class sizes are given.

The binary difference \eqref{binary_difference_formula}, plotted on log-scale against the iteration count, is depicted for each experiment in Figure \ref{fig:binary_difference}. For experiments without any size information, the average binary difference
tends to less than $10^{-16}$, which is vanishingly low and more or less indicates that an exact global minimizer has been obtained.
For experiments with size constraints or penalty terms, the binary difference also gets very low, although not as low. This
indicates convergence to at least a very close approximation of the global minimizer. These observations agree
well with the theoretical results in Section \ref{analysis dual}, where the strongest results were also obtained in case of no size information.

Note that a lot more iterations than necessary have been used in the binary difference plots. In practice, the algorithm reaches sufficient stability in 100-300 iterations. The CPU times, summerized in Table \ref{table:timing}, indicate a fast convergence of the new algorithm, much faster than GL, although not quite as fast as the MBO scheme. It must be noted that MBO is an extremely fast front propagation algorithm that only uses a few (e.g. 10) iterations, but its accuracy is limited due to the large step sizes. A deeper discussion on the number of iterations needed to reach the exact solution after thresholding will be given at the end of the next section on point cloud segmentation.

\subsection{Segmentation of 3D point clouds}\label{Segmentation of point clouds}

The energy function \eqref{original problem with fidelity} that combines region homogeneity terms and dissimilarity across region boundaries will be demonstrated for segmentation of unstructured 3D point clouds, where each point is a vertex in $V$. Point clouds arise for instance through laser-based range imaging or multiple view scene reconstruction. The results of point cloud segmentation are easy to visualize and the choice of each term in the energy function will have a clear intuitive meaning that may be translated to other graph-based classification problems in the future. We focus especially on point clouds acquired through the concept of laser detection and ranging (LaDAR) in outdoor scenarios. A fundamental computer vision task is to segment such scenes into classes of similar objects. Roughly, some of the most common object classes in an outdoor scene are the ground plane, vegetation and human-made "objects" with a certain regular structure.


\subsubsection{Construction of the energy function}

We construct the graph by connecting each node to its k nearest neighbors (kNN) based on the Euclidian distance as described at the beginning of Section \ref{Segmentation as optimization problem over graph}. In experiments, we set $k=20$. We construct
region terms that favor homogeneity of geometrical features based on a combination of point coordinates, normal vectors and variation of normal vectors. The construction is a concrete realization of the general region terms introduced in \cite{lezoray,LEL14,TLE15}. We also propose to use a contour term that favors alignment of the boundaries of the regions at "edges", indicated by sharp discontinuities of the normal vectors. Our model can be seen as a point cloud analogue of variational models for traditional image segmentation, that combine region and edge based features in a single energy functional \cite{bresson2007fast,gilboa2,DBLP:journals/siamis/JungPC12}. In contrast to the work \cite{DBLP:conf/cvpr/AnguelovTCKGHN05,MVH08,Triebel06robust3d} our model does not rely on training data.




Normal vectors in a point cloud can be estimated from principal component analysis locally around each point, as in e.g. \cite{Di12,LLZ13,LEL14}. For each point $x \in V$, let $y^1,...,y^m$ denote the set of neighboring points and define for notational convenience $y^0 = x$. Define the normalized vectors
$\bar{y}^i = y^i - \text{mean}(y^0,y^1,...,y^m)$ for $i=0,1,...,m$ and construct the matrix
\begin{equation}\label{correlation matrix}
\textbf{Y} = [\bar{y}^0 \bar{y}^1 \bar{y}^2 ... \bar{y}^m].
\end{equation}
Let $\textbf{v}^1(x),\textbf{v}^2(x),\textbf{v}^3(x)$ be the eigenvectors and \newline $\lambda^1(x),\lambda^2(x),\lambda^3(x)$ be
the eigenvalues of the correlation matrix $\textbf{Y} \textbf{Y}^T$. The first eigenvector $\textbf{v}^1(x)$ points in the direction of least variation between the points $\bar{y}^1 , ... , \bar{y}^m$
and the first eigenvalue $\lambda^1(x)$ indicates the variation along the direction of $\textbf{v}^1(x)$.

The variable $\textbf{v}^1(x)$ is consequently a discrete estimation of the
normal vector at $x$ and 
the first eigenvalue $\lambda^1(x)$ indicates to which extend the normal vectors vary locally around the point $x$. If all the points were laying on a plane, then $\lambda^1(x)$ would be zero and $\textbf{v}^1(x)$ would be the normal vector of the
plane. 

The region term for region $V_i$ can be constructed to be small at the point $x$ if the value of $\lambda^1(x)$ is close to the expected value $\lambda_i$ of region $i$, and be large otherwise. This can be achieved by requiring the following term to be small
\begin{equation}\label{orientation variation fidelity term}
 \big| \lambda^1(x) - \lambda_i \big|^2 \, , \quad \forall x \in V, \quad i = \{v, h, g\}.
\end{equation}
For instance, $\lambda^1_v$, $\lambda^1_h$ and $\lambda^1_g$ for vegetation, human-made objects and the ground plane can be estimated from measurements. Note that their particular values depend on characteristics of the LaDAR, such as the angular scanning resolution, depth resolution etc. If an estimate of $\lambda_i$ is not known, $\lambda_i$ could be part of the minimization problem in a similar way to the mean intensity values in the Chan-Vese model \cite{CV01}.

Furthermore, the region terms can be constructed for discriminating regions where the normal vector are oriented either parallel with or perpendicular to a specific direction $\textbf{n}^i$ by
 requiring the following terms to be small, respectively
$$- |\textbf{v}^1(x) \cdot \textbf{n}^i|, \quad \quad |\textbf{v}^1(x) \cdot \textbf{n}^i|.$$
For instance, the normal vectors of the ground plane are expected to point predominantly in the upward direction.
The ground plane can also be
characterized by its height, defined by the $x_3$-coordinate of the points, which is generally lower than the heights of other objects in the nearby surroundings. Assuming a rough estimate of the local height of the ground plane $h^*(x)$ at the point $x$ is known, the fidelity term \eqref{original problem with fidelity} can be modified to take into account both normals vectors and height by requiring the following term to be small
\begin{equation}\label{fidelity nj}
 - |\textbf{v}^1(x) \cdot \textbf{n}^i(x)| + H(x_3, h^*(x)),
\end{equation}
where $H$ is an increasing function in $x_3$ and, in addition, $H(h^*(x),h^*(x)) = 0$. We have used the term $H(x_3, h^*(x)) = \theta (x_3 - h^*(x))$ and simply estimated $h^*(x)$ as the average $x_3$ coordinate of the points in the neighborhood of $x$.

\begin{figure}
  \begin{tabular}[t]{c}
   \includegraphics[width=0.45\textwidth]{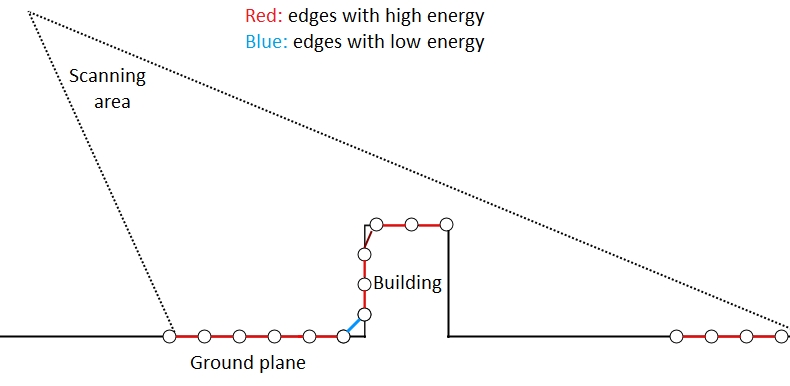}
  \end{tabular}
 \caption{\label{figure ground plane illustration}Illustration of construction of a graph to separate the ground plane from human-made structures, view point from the side. The edges are assigned a low energy at convex parts of the scene, marked in light blue, making it favorable to place the boundary between the regions at such locations.}
\end{figure}

\begin{figure*}
 \begin{tabular}[b!]{c}
    \subfigure[Primal energies]
   {\mbox{
   \includegraphics[width=0.33\textwidth]{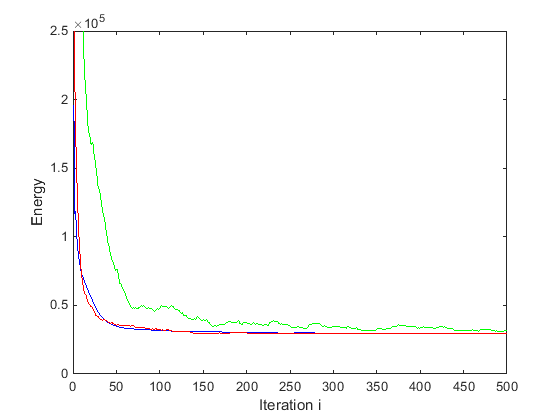}}}
       \subfigure[Log of $\frac{||E^i - E^*_{relaxed}||}{E^*_{relaxed}}$]
   {\mbox{
   \includegraphics[width=0.33\textwidth]{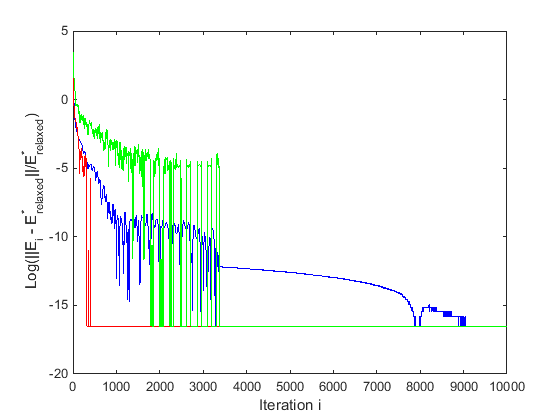}}}
   \subfigure[Binary difference \eqref{binary_difference_formula}]
   {\mbox{
   \includegraphics[width=0.33\textwidth]{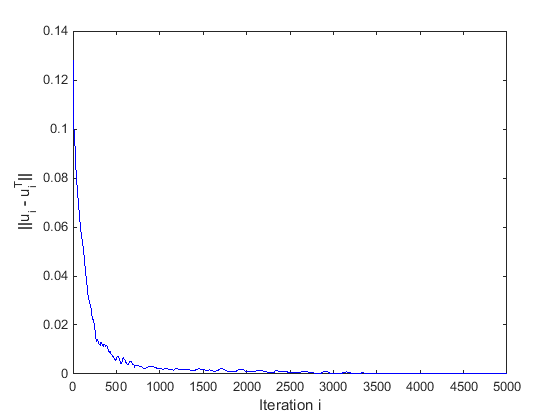}}}
   \end{tabular}
 \caption{\label{fig:energy convergence}(a)-(b) Energy evolution of $u$ (blue), $u^T$ with threshold scheme \eqref{rounding} (red), and $u^T$ with threshold scheme \eqref{binary primal solution} (green) for the experiment in Figure \ref{figure slangeblokk}.}
\end{figure*}



The weight function $w$ is constructed to encourage spatial grouping of the points, and so that it is favorable to align the border between regions at locations where the normal vectors change from pointing upwards to pointing outwards, i.e. where the scene is convex-shaped. On the contrary, locations where the scene is concave, such as the transition from the side of buildings to the roof, should be unfavorable for the region boundaries. Such assumptions can be incorporated by modifying the Gaussian weight function \eqref{w1} as follows:
 \begin{equation}
            w(x,y)= e^{-\frac{d(x,y)^{2}}{\sigma^{2}} + \gamma \frac{v^1_3(y)-v^1_3(x)}{d(x,y)}\text{SIGN}(y_1 - x_1)}
\label{w2}
\end{equation}
Here $v_1^1(x)$ and $v_3^1(x)$ are the first and third components of the vector $\textbf{v}^1(x)$, and a coordinate system has been assumed where the positive $x_1$ axis points outwards from the view direction and the positive $x_3$ axis points upwards. An illustration is given in Figure \ref{figure ground plane illustration}, where edges at convex parts of the scene are given a low energy value, marked by the color code of light blue.


Taking the above information into account, an example of how the different region terms can be constructed for the ground plane, human-made structures and vegetation, respectively, are
\begin{align}
f_g(x) & = (1-C)\big| \lambda^1(x) - \lambda_g \big|^2 \nonumber \\
 + C & \big\{ - |\textbf{v}^1(x) \cdot \textbf{n}^g(x)| + H(x_2,h^*(x)) \big\}. \label{fidelity nj1} \\
f_h(x) & = (1-C) \big| \lambda^1(x) - \lambda_h \big|^2 \nonumber \\
& + C |\textbf{v}^1(x) \cdot \textbf{n}^g(x)| , \label{fidelity nj2} \\
f_v(x) & = C \big| \lambda^1(x) - \lambda_v \big|^2. \label{fidelity nj3}
\end{align}
Here, $C \in (0,1)$ is a parameter that balances considerations between variation and direction/height. In experiments, we set
$\lambda_g = \lambda_h$ and set $C$ to a low value so that only vegetation is distinguished from other regions by the value of $\lambda_1$.
In some experiments, we also use two regions for vegetation with two different values of $\lambda_i$. This makes it possible to distinguish different kinds of vegetation, those with leaves or needles tend to have a lower mean value $\lambda_i$ than those without them. 
Smoke can also characterized by its irregular surface, and its region term constructed as  \eqref{orientation variation fidelity term} with an intermediate value of $\lambda_i$.

\begin{figure}
  \begin{tabular}[b]{c}
  \subfigure[Scanning area]
   {\mbox{
   \includegraphics[width=0.45\textwidth]{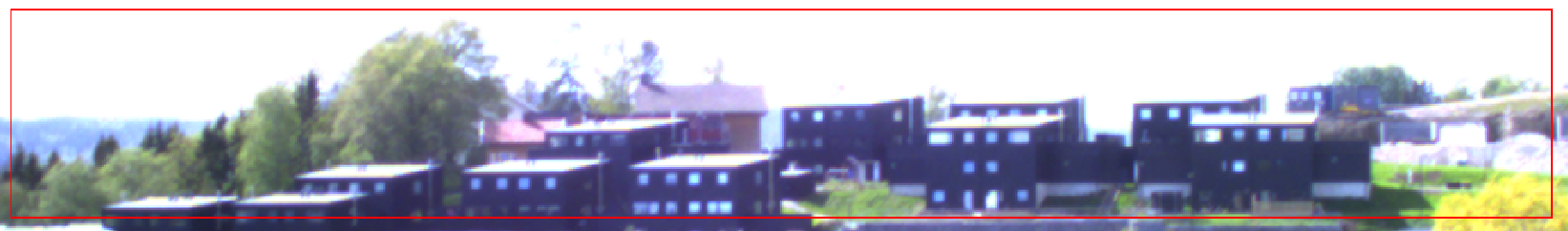}}}\\
     \subfigure[Segmentation, view from front]
   {\mbox{
   \includegraphics[width=0.45\textwidth]{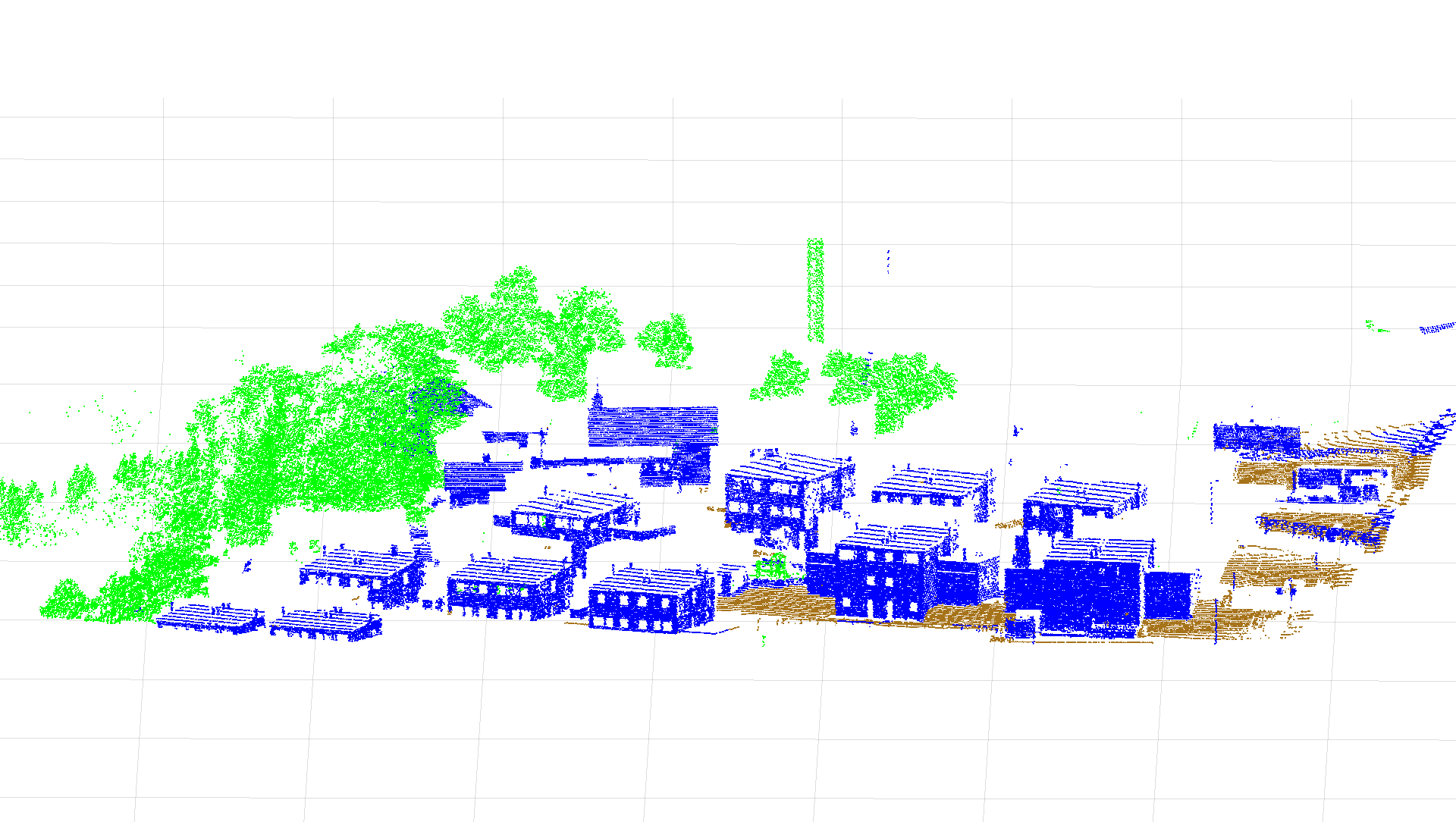}}}\\
     \subfigure[Segmentation, view from top]
   {\mbox{
   \includegraphics[width=0.45\textwidth]{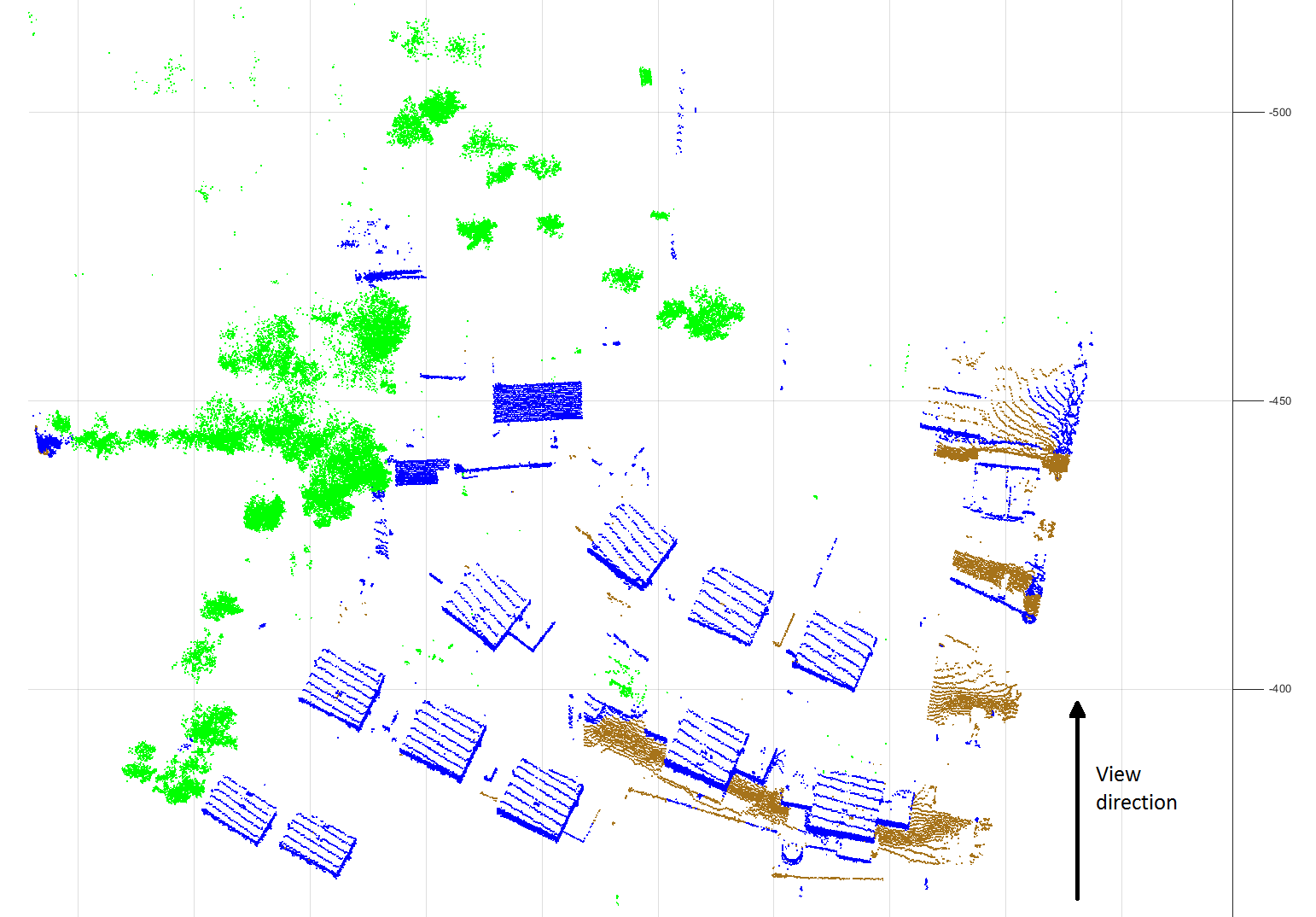}}}
  \end{tabular}
  \caption{\label{figure rekkehus}(a) Scanning area of LaDAR. (b)-(c) Segmentation of acquired point cloud, consisting of $93641$ points, into 3 regions: ground plane (brown), vegetation (green) and human-made
 objects (blue).}
\end{figure}
\begin{figure}
  \begin{tabular}[t]{c}
    \subfigure[Scanning area]
   {\mbox{
   \includegraphics[width=0.45\textwidth]{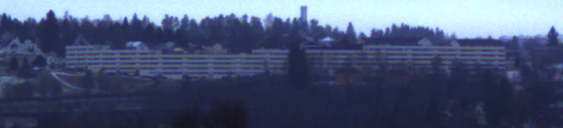}}}\\
        \subfigure[Segmentation, view from front]
   {\mbox{
   \includegraphics[width=0.45\textwidth]{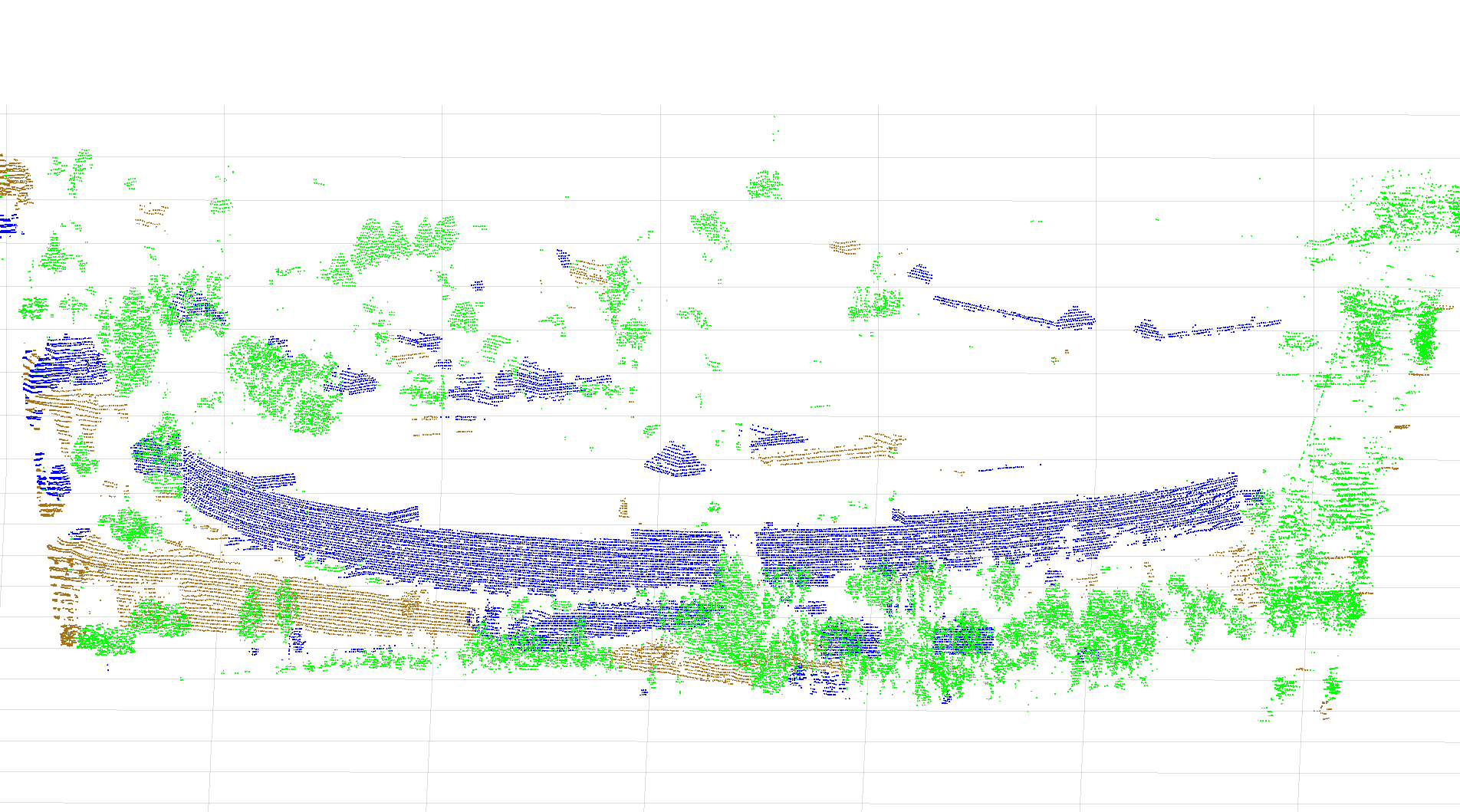}}}\\
        \subfigure[Segmentation, view from top]
   {\mbox{
   \includegraphics[width=0.45\textwidth]{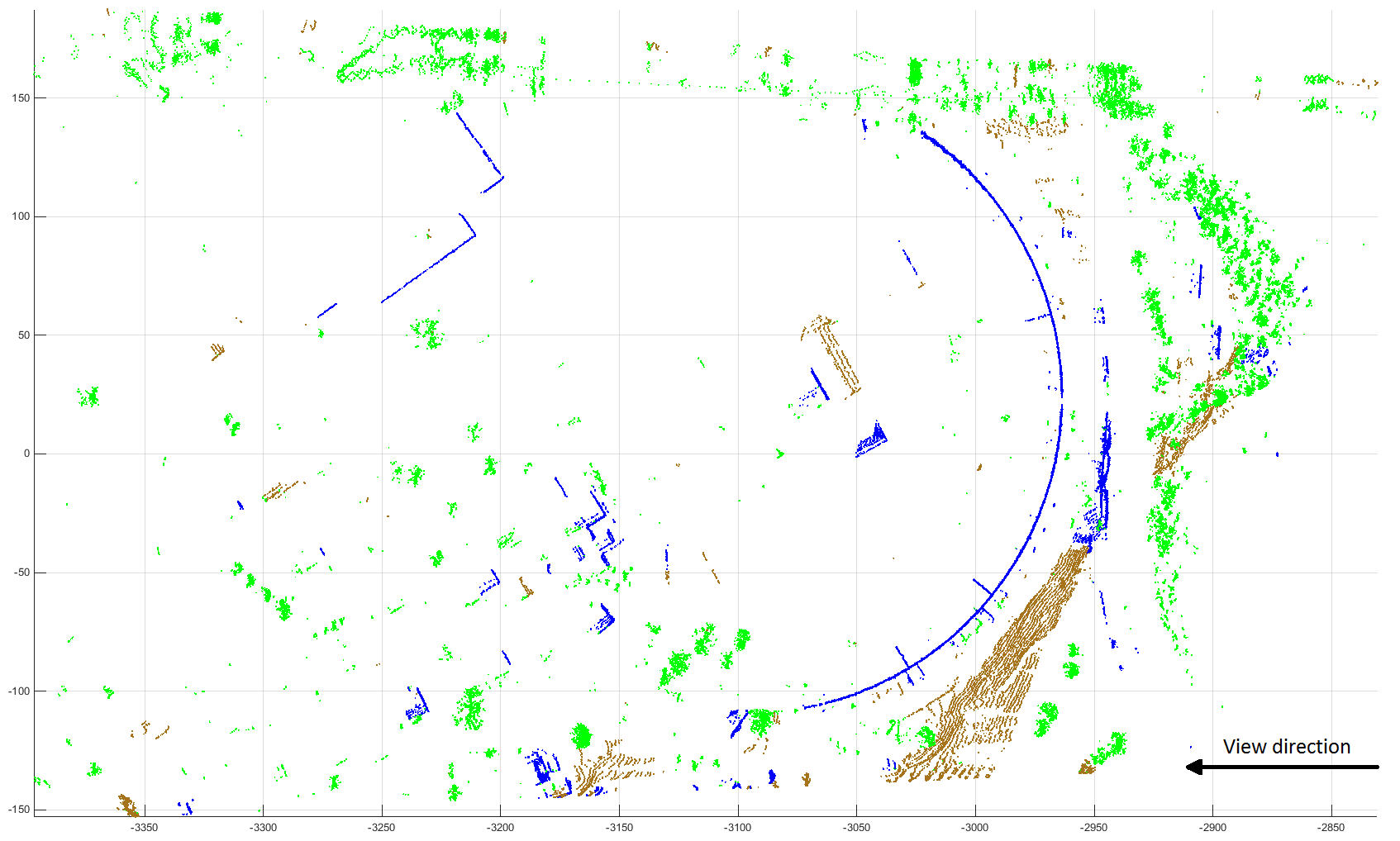}}}
  \end{tabular}
 \caption{\label{figure slangeblokk}
 (a) Scanning area of LaDAR. (b)-(c) Segmentation of acquired point cloud, consisting of $80937$ points, into 3 regions: ground plane (brown), vegetation (green) and human-made objects (blue).}
\end{figure}

\subsubsection{Experiments}

\begin{figure}
 \begin{tabular}[t!]{c}
   \includegraphics[width=0.45\textwidth]{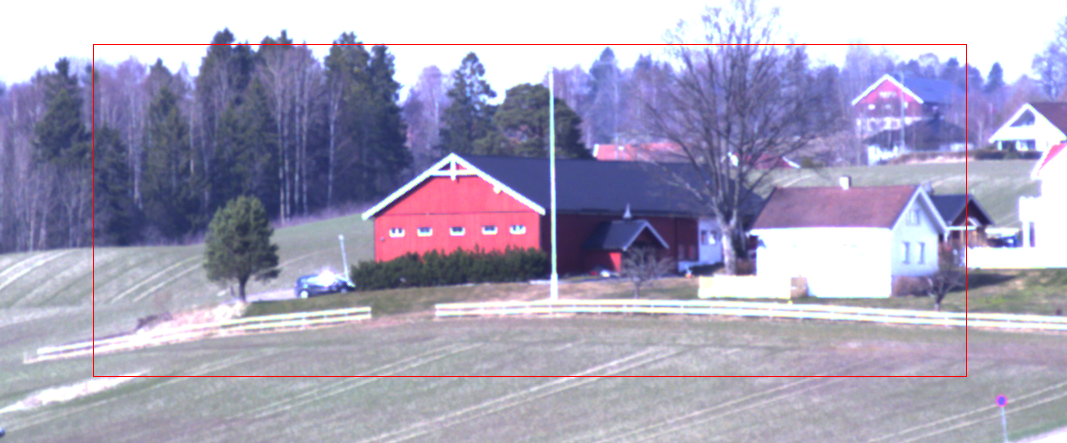}\\
   \includegraphics[width=0.45\textwidth]{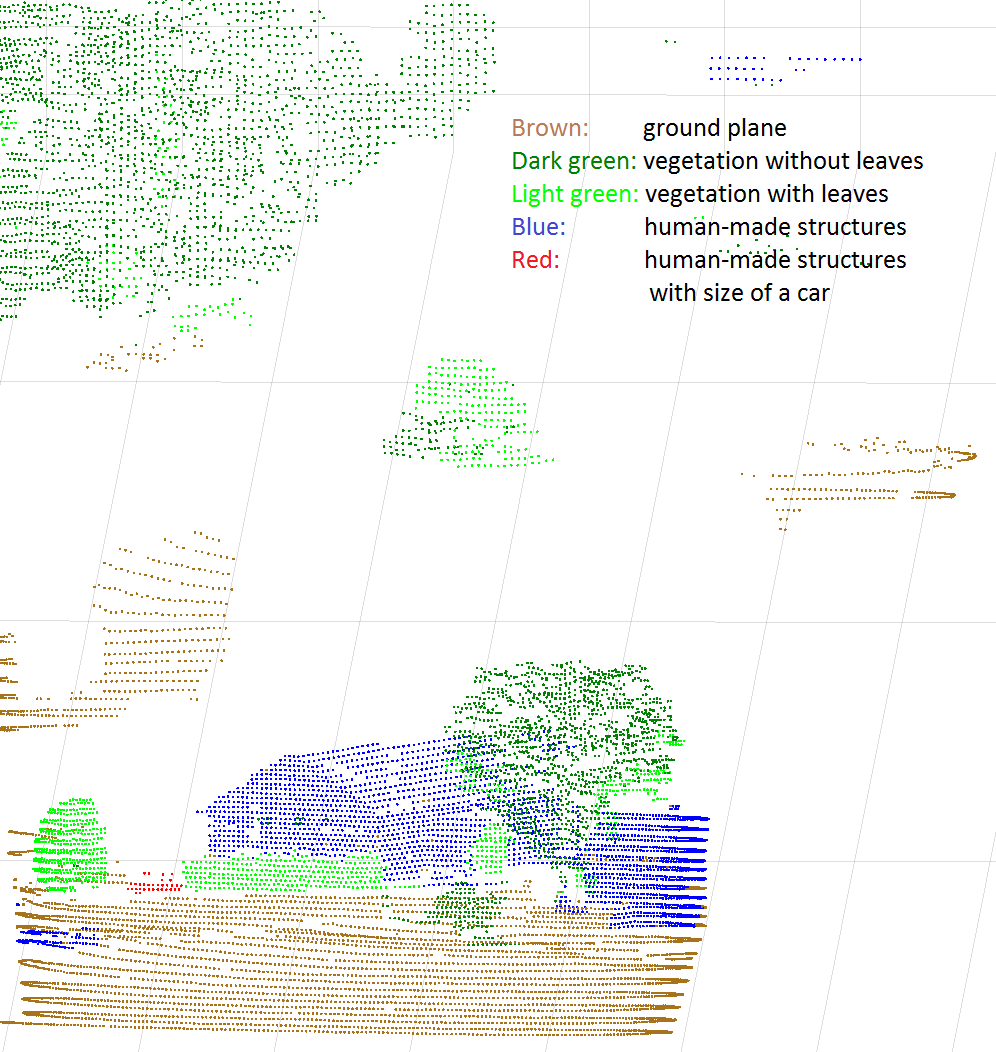}
     \end{tabular}
 \caption{\label{fig:farm}
 Top: Scanning area of LaDAR. Bottom: Segmentation of acquired point cloud, consisting of $14002$ points, into 4 regions: ground plane (brown), and human-made
 objects (blue), vegetation with (light green) and without (dark green) leaves/needles. }
\end{figure}

Some illustrative experiments are shown in Figures \ref{figure rekkehus} and \ref{figure slangeblokk}. Ordinary photographs of the scenes are shown on the top and the red rectangles indicate the areas that have been scanned by the LaDAR. The point clouds have been segmented into three regions as described above and the results are visualized by brown color for points assigned to the ground plane region, green color for points assigned to the vegetation region and blue color for points assigned to the region of human-made objects. In Figure \ref{fig:farm}, vegetation with and without leaves are indicated by dark and light green respectively.

\begin{figure}
 \begin{center}
 \begin{tabular}[t!]{c}
    \subfigure[Scanning area]
   {\mbox{
   \includegraphics[width=0.45\textwidth]{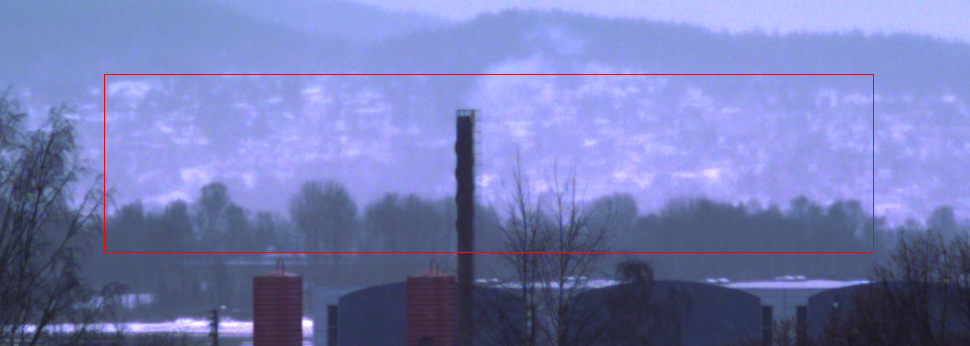}}}\\
    \subfigure[Segmentation result]
   {\mbox{
   \includegraphics[width=0.45\textwidth]{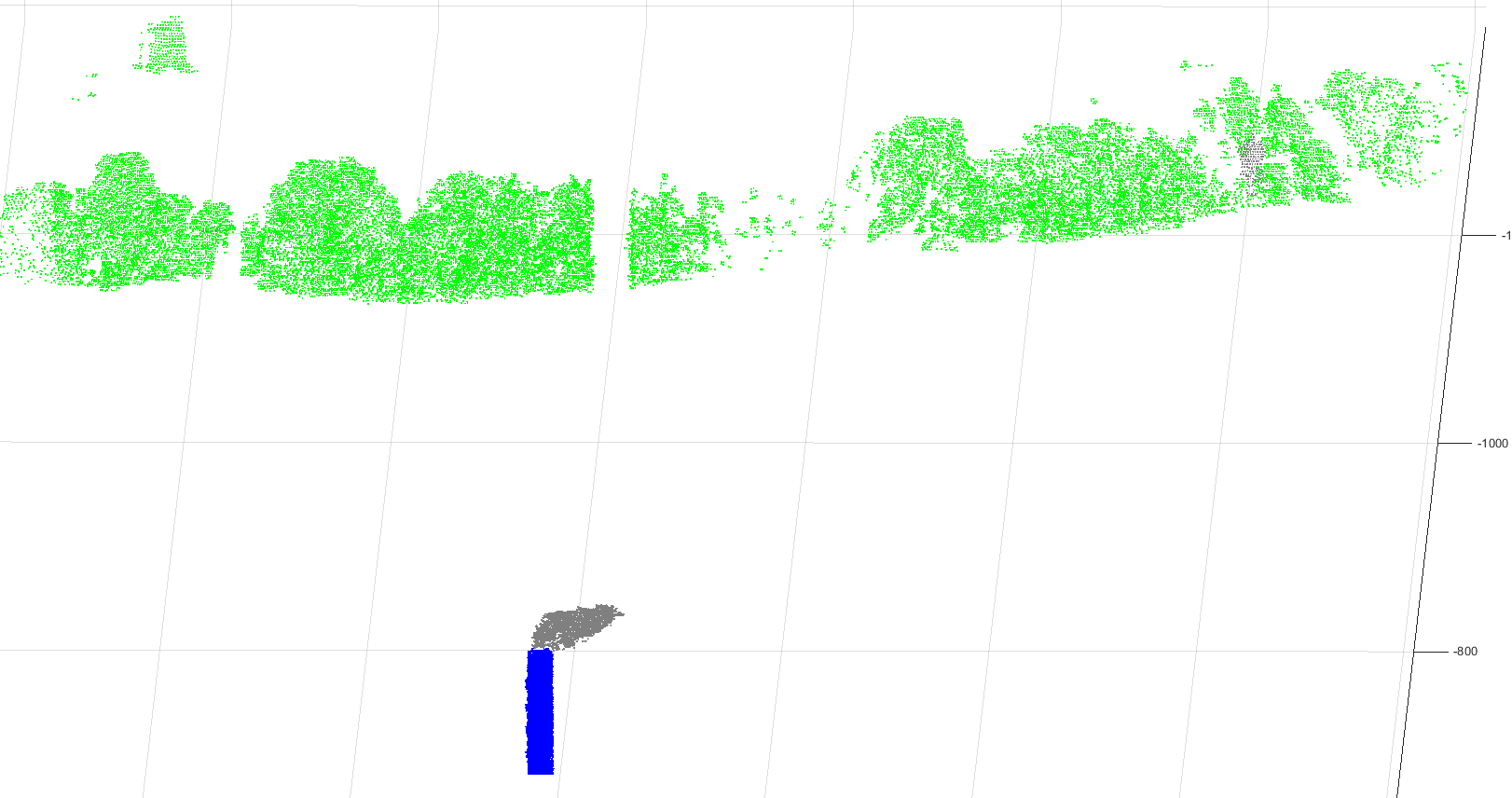}}}\\
     \end{tabular}
 \end{center}
 \caption{\label{fig:pipe}Top: Scanning area of LaDAR. Bottom: Segmentation of a point cloud (81551 points) into smoke (gray), vegetation (green) and human-made structures (blue).}
\end{figure}
%

It can be observed that the algorithm leads to consistent results even though these scenes are particularly challenging because the tilt and height of the ground plane vary highly over the scene due to the hilly landscape, and some of the trees and bushes are completely aligned with and touches the buildings.
Note that buildings hidden behind vegetation get detected since the laser pulses are able to partially penetrate through the leaves. 
A misassignment can be observed in the middle of Figure \ref{figure slangeblokk}, where only the roof of one of the buildings is visible due to occlusions. Since no points are observed from the wall of the building, the roof gets assigned to the ground plane region. Some large rocks on figure \ref{figure slangeblokk} also get assigned to the blue region due to their steep and smooth surfaces.

%
%

As was the case for experiments involving semi-supervised classification, the approximation errors of the convex relaxation practically vanish. Figure \ref{fig:energy convergence}(c) depicts the binary difference \eqref{binary_difference_formula} as a function of the iteration count in the experiment shown in Figure \ref{figure slangeblokk}. As can be seen, the solution of the convex relaxation converges to a binary function; after 10000 iterations, the average binary difference \eqref{binary_difference_formula} was $5.74*10^{-10}$. Note, however, that a lot less iterations are necessary before the thresholded function stabilizes at the global minimum. Figure \ref{fig:energy convergence} (left) depicts the energy evolution as a function of the iteration count for the relaxed solution (blue), thresholded solution with scheme \eqref{rounding} (red) and with scheme \eqref{binary primal solution} (green). Figure \ref{fig:energy convergence} (right) depicts a log plot of the absolute energy precision
$\frac{||E^i - E^*_{relaxed}||}{E^*_{relaxed}}||$, where $E^*_{relaxed}$ is the global minimum of the relaxed problem, estimated by 10000 iterations of the algorithm. $E^i$ is the energy at iteration $i$ of the relaxed solution (blue), thresholded solution with scheme \eqref{rounding} (red) and thresholded solution with scheme \eqref{binary primal solution} (green). This plot demonstrates that the binary solution obtained by the thresholding scheme \eqref{rounding} stabilizes after about 300 iterations, after which the energy is within $10^{-16}$ of the energy of the ground truth solution of the relaxed problem estimated at iteration 10000. The threshold scheme \eqref{binary primal solution} takes more iterations before stabilizing, but also eventually converges to the correct solution after about 3500 iterations. The CPU times of the experiments were in the range 5-15 seconds on an Intel i5-4570 3.2 Ghz CPU for point clouds with around 80000 points. For comparison, the inference step of the related MRF approaches \cite{DBLP:conf/cvpr/AnguelovTCKGHN05,Triebel06robust3d,MVH08} took around 9 minutes for a scan with around 30000 points, as reported in \cite{MVH08}, but of course on older hardware. The proposed algorithm is also suitable for parallel implementation on GPU as discussed at the end of Section \ref{sec algorithm}.

\section{Conclusions}

Variational models on graphs have shown to be highly competitive for various data classification problems, but are inherently difficult to handle from an optimization perspective, due to \newline NP-hardness except in some restricted special cases. This work has developed an efficient convex algorithmic framework for a set of classification problems with multiple classes involving graph total variation, region homogeneity terms, supervised information and certain constraints or penalty terms acting on the class sizes. 
Particular problems that could be handled as special cases included semi-supervised classification of high-dimensional data and unsupervised segmentation of unstructured 3D point clouds. The latter involved minimization of a novel energy function enforcing homogeneity of point coordinate based features within each region, together with a term aligning the region boundaries along edges. Theoretical and experimental analysis revealed that the convex algorithms were able to produce vanishingly close approximations to the global minimizers of the original problems in practice.

Experiments on benchmark datasets for semi-supervised classification resulted in higher accuracies of the new algorithm compared to related local minimization approaches. The accuracies were also highly competitive against a wide range of other established algorithms. The advantages of the proposed algorithm were particularly prominent in case of sparse or non-uniformly distributed training data. The accuracies could be improved further if an estimate of the approximate class sizes were given in advance. Experiments also demonstrated that 3D point clouds acquired by a LaDAR in outdoor scenes could be segmented into object classes with a high degree of accuracy, purely based on the geometry of the points and without relying on training data. The computational efficiency was at least an order of magnitude faster than related work reported in the literature.

In the future, it would be interesting to investigate region homogeneity terms for general unsupervised classification problems. In addition to avoiding the problem of trivial global minimizers, the region terms may improve the accuracy compared to models based primarily on boundary terms. Region homogeneity may for instance be defined in terms of the eigendecomposition of the covariance matrix or graph Laplacian.

\appendix

\section{Proof of Theorem \ref{theorem two components}}\label{appendix 1}

To aid the proof of Theorem \ref{theorem two components}, we first give the following lemma, which is a graph extension of Proposition 4 given in \cite{BYT2011} for image domains.
\begin{lemma}\label{lemma threshold}
Assume that for a function $u\,: V \mapsto [0,1]$, $q^*$ maximizes
$$q^* = \argmax_{\norm{q}_{\mathcal{E},\infty} \leq 1} \sum_{x \in V} u(x) (\diver_w q)(x)$$
Define the thresholded function
\begin{equation}\label{eqn: u}
u^\alpha(x) \; = \; \left\{ \begin{array}{ll}
                           1 & \textrm{ if $u(x) \geq \alpha$} \, \\
                           0 & \textrm{ otherwise}
                    \end{array}
\right. .
\end{equation}
For almost any threshold level $\alpha \in (0,1]$, $q^*$ also maximizes
$$q^* = \argmax_{\norm{q}_{\mathcal{E},\infty} \leq 1} \sum_{x \in V} u^\alpha(x) (\diver_w q)(x)$$
\end{lemma}
\begin{proof}
The coarea formula on graphs says that
$$\sum_{x \in V} |\nabla_w u(x)| = \int_{0}^1 \sum_{x \in V} |\nabla_w u^\alpha(x)| \, d \alpha,$$
see for instance appendix B of \cite{GGOB13} for a proof. Together with the fact that $u(x) = \int_0^{u(x)} d \alpha = \int_0^{1} u^\alpha(x)d \alpha$, we can deduce that
\begin{align*}
& \int_{0}^1 \sum_{x \in V} u^\alpha(x) (\diver_w q^*)(x) \, d \alpha \\
= &\sum_{x \in V} \big(\int_{0}^1  u^\alpha(x) \, d \alpha \big) (\diver q^*)(x) \\
= &\sum_{x \in V} u(x) (\diver_w q^*)(x)\\
 = &\sum_{x \in V} |\nabla_w u(x)| = \int_{0}^1 \sum_{x \in V} |\nabla_w u^\alpha(x)| \, d \alpha\\
= & \int_{0}^1 \big( \sup_{\norm{q}_{\mathcal{E},\infty} \leq 1} \sum_{x \in V} u^\alpha(x) (\diver q)(x) \, d \alpha \big).
\end{align*}
Since in general
\begin{align*}
\sup_{\norm{q}_{\mathcal{E},\infty} \leq 1}& \sum_{x \in V} u^\alpha(x) (\diver q)(x) \, d \alpha \\
\geq & \sum_{x \in V} u^\alpha(x) (\diver q^*)(x) \, d \alpha,
\end{align*}
the above equality can only be true provided that
\begin{align*}
\sup_{\norm{q}_{\mathcal{E},\infty} \leq 1} & \sum_{x \in V} u^\alpha(x) (\diver q)(x) \\
= & \sum_{x \in V} u^\alpha(x) (\diver q^*)(x) ,
\end{align*}
for almost every $\alpha \in (0,1]$.
\end{proof}
Utilizing Lemma \ref{lemma threshold}, we will now prove Theorem \ref{theorem two components}:
\begin{proof}
By the assumptions of the theorem, for a finite number of connected components in the graph, the minimizer $I_\text{min}(x)$ contains two indices. Assume without loss of generality that $V_{k,j} \subset V$ is one such connected component
where $I_m(x) = k,j$ for all $x \in V_{k,j}$. That is, for any two nodes $x$ and $y$ in $V_{k,j}$, there is a path of edges
$(x,z_1),(z_1,z_2),...,(z_n,y) \subset E$ such that $z_1,...,z_n \in V_{k,j}$.

Let $u^*$ be any primal solution for which $(u^*;q^*)$ is a
primal-dual pair. By Theorem \ref{theo:consist_unq}, $u^*$ must in $V_{k,j}$ satisfy
\begin{equation}\label{relation two minimums}
u^*_k(x) + u^*_j(x) = 1, \quad u^*_i(x) = 0, \quad \text{for} \; i \neq k,j,
\end{equation}
For an arbitrary threshold level $\alpha \in (0,1)$ construct the binary function
\begin{equation}\label{eqn: u}
u^\alpha_k(x) \; = \; \left\{ \begin{array}{ll}
                           1 & \textrm{ if $u^*_k(x) \geq \alpha$} \, \\
                           0 & \textrm{ otherwise}
                    \end{array}
\right. .
\end{equation}
From \eqref{relation two minimums}, we can write $u^*_j(x) = 1 - u^*_k(x)$ in $V_{k,j}$,
and together with \eqref{eqn: u} it follows that $1 - u^\alpha_k(x) = u^{1-\alpha}_j(x)$ in $V_{k,j}$.

Construct now the function $u^t \, : \, V \mapsto \mathbb{R}^n$ as follows:
\begin{equation}\label{eqn: ut}
u^t(x) = u^*(x)\; \text{for}\; x \in V \backslash V_{k,j}
\end{equation}
\begin{equation}\label{eqn: ut2}
u^t_i(x) \; = \; \left\{ \begin{array}{ll}
                           u^\alpha_k(x) & \textrm{if $i=k$} \, \\
                           u^{1-\alpha}_j(x) & \textrm{if $i=j$}\\
                           0 & \textrm{if $i \neq k,j$}
                    \end{array}
\right. \; \text{for}\; x \in V_{k,j}
\end{equation}
For the given $q^*$, we have that
\begin{align}
& E(u^*,q^*) \nonumber \\
= & \sum_{i\in I} \sum_{x \in V}  u_i^*(x) \big\{C_i(x) \, +\, (\diver_w q_i^*)(x) \big\} \nonumber \\
= & \sum_{i\in I \backslash \{k,j\}} \sum_{x \in V \backslash V_{k,j}}  u_i^*(x) \big\{C_i(x) \, +\, (\diver_w q_i^*)(x) \big\}\nonumber \\
+ & \sum_{x \in V_{k,j}} u_k^*(x) \big\{C_k(x) \, +\, (\diver_w q_k^*)(x) \big\} \nonumber \\
+ & \sum_{x \in V_{k,j}}  u_j^*(x) \big\{C_j(x) \, +\, (\diver_w q_j^*)(x) \big\}\nonumber \\
= & \sum_{i\in I \backslash \{k,j\}} \sum_{x \in V \backslash V_{k,j}}  u_i^*(x) \big\{C_i(x) \, +\, (\diver_w q_i^*)(x) \big\}\nonumber \\
+ & \sum_{x \in V_{k,j}} \big( u_k^*(x) + (1-u_k^*(x)) \big) \big\{C_k(x) \, +\, (\diver_w q_k^*)(x) \big\} \nonumber \\
= & \sum_{i\in I \backslash \{k,j\}} \sum_{x \in V \backslash V_{k,j}}  u_i^*(x) \big\{C_i(x) \, +\, (\diver_w q_i^*)(x) \big\}\nonumber \\
+ & \sum_{x \in V_{k,j}} \big( u_k^\alpha(x)+ (1-u_k^\alpha(x))\big) \big\{C_k(x) \, +\, (\diver_w q_k^*)(x) \big\}\nonumber
\end{align}
\begin{align}
= & \sum_{i\in I \backslash \{k,j\}} \sum_{x \in V \backslash V_{k,j}}  u_i^*(x) \big\{C_i(x) \, +\, (\diver_w q_i^*)(x) \big\}\nonumber\\
+ & \sum_{x \in V_{k,j}}  u^\alpha_k(x) \big\{C_k(x) \, +\, (\diver_w q_k^*)(x) \big\},\nonumber \\
+ & \sum_{x \in V_{k,j}}  u^{1-\alpha}_j(x) \big\{C_j(x) \, +\, (\diver_w q_j^*)(x) \big\},\nonumber \\
= & \, E(u^t,q^*)
\end{align}
where we have used that $C_k(x) \, +\, (\diver_w q_k^*)(x) = C_j(x) \, +\, (\diver_w q_j^*)(x)$ in
$V_{k,j}$. 
By applying Lemma \ref{lemma threshold} on the last two terms with threshold level $\alpha$ and $1-\alpha$ respectively, it can be deduced that
$$\sup_{q \in S^n_\infty}E(u^t,q) = \sup_{q \in S^n_\infty} E(u^*,q) = E(u^*,q^*) = E(u^t,q^*) $$
Consequently $(u^t,q^*)$ is an optimal primal-dual pair.

Assume now there is another connected component $V^2_{k_2,j_2} \subset V$ where $I_\text{min} = \{k_2,j_2\}$. By setting $u^* = u^t$ and repeating all arguments above, it follows that $u^*$ can be thresholded in $V^2_{k_2,j_2}$ to yield a binary minimizer in $V^2_{k_2,j_2}$. The same process
can be repeated for all connected components until a binary minimizer is obtained over the whole domain $V$. By Proposition
\ref{primal: exact global minimizer}, such a binary function is a global minimizer of the original non-convex problem.

\end{proof}

\bibliographystyle{plain}
\bibliography{References_short}



\end{document}